\documentclass[a4paper,12pt]{article}
\title{
Quatre applications du lemme de Zalcman\\
 \`a la dynamique complexe
}
\author{
Tomoki Kawahira
\thanks{
financ\'e en partie par la Fondation Sumitomo, 
la Fondation Chubei Itoh, 
et la JSPS}
\\ {Universit\'e de Nagoya}
}

\usepackage{enumerate}
\usepackage[dvipdfmx,pdftex]{graphicx}
\usepackage{amsmath}
\usepackage{amscd}
\usepackage{amssymb}
\usepackage{amsfonts}
\usepackage{theorem}
\usepackage[english,french]{babel}

\newtheorem{thm}{Th\'eor\`eme}[section]
\newtheorem{prop}[thm]{Proposition}
\newtheorem{lem}[thm]{Lemme}
\newtheorem{cor}[thm]{Corollaire}
\theorembodyfont{\rmfamily}
\newtheorem{pf}{D\'emonstration.}

\newtheorem{df}{D\'efinition.}

\newcommand{\thmref}[1]{Th\'eor\`eme \ref{#1}}
\newcommand{\propref}[1]{Proposition \ref{#1}}

\newcommand{\lemref}[1]{Lemme \ref{#1}}
\newcommand{\figref}[1]{Figure \ref{#1}}

\newcommand{\parag}[1]{
\par\medskip
\noindent {\bf #1}
}

\newcommand{\C}{{\mathbb{C}}}
\newcommand{\B}{{\mathbb{B}}}
\newcommand{\Cstar}{{\mathbb{C}^\ast}}
\newcommand{\Chat}{{\widehat{\mathbb{C}}}}

\newcommand{\D}{{\mathbb{D}}}

\newcommand{\N}{{\mathbb{N}}}
\newcommand{\abs}[1]{{\left| #1 \right|}}
\newcommand{\paren}[1]{{\left( #1 \right)}}
\newcommand{\braces}[1]{{\left\{ #1 \right\}}}
\newcommand{\gauss}[1]{{\left [ #1 \right ]}}

\newcommand{\dist}{\mathrm{dist}\,}

\newcommand{\al}{{\alpha}}
\newcommand{\gam}{{\gamma}}
\newcommand{\lam}{{\lambda}}
\newcommand{\del}{{\partial}}

\newcommand{\e}{\epsilon}
\newcommand{\cc}{\circ}
\newcommand{\cA}{{\mathcal{A}}}
\newcommand{\cF}{{\mathcal{F}}}
\newcommand{\cJ}{{\mathcal{J}}}
\newcommand{\cK}{{\mathcal{K}}}
\newcommand{\cM}{{\mathcal{M}}}
\newcommand{\cR}{{\mathcal{R}}}
\newcommand{\cL}{{\mathcal{L}}}
\newcommand{\cU}{{\mathcal{U}}}
\newcommand{\cV}{{\mathcal{V}}}
\newcommand{\cZ}{{\mathcal{Z}}}
\newcommand{\st}{\,:\,}
\newcommand{\QED}{\hfill $\blacksquare$}
\newcommand{\zhat}{\hat{z}}
\newcommand{\cLM}{\cL\cM}
\newcommand{\psihat}{\hat{\psi}}
\newcommand{\phihat}{\hat{\phi}}
\newcommand{\cZA}{\cA^{\cZ}}
\newcommand{\cZG}{\mathcal{G^Z}}
\newcommand{\fhat}{\widehat{f}}
\newcommand{\Lam}{\varLambda}
\newcommand{\Aff}{\mathrm{Aff}}

\begin{document}

\maketitle

\begin{abstract}
Nous donnons quatre applications du lemme de Zalcman
\`a la dynamique des fractions rationnelles sur la sph\`ere de Riemann: 
un analogue param\'etrique de la d\'emonstration de la densit\'e des cycles r\'epulsifs; 
la ressemblance de l'ensemble de Mandelbrot avec les ensembles de Julia; 
une construction de la lamination de Lyubich-Minsky et d'une variante;
et une caract\'erisation unifi\'ee
des points coniques de Lyubich-Minsky et ceux de Martin-Mayer. 
\begin{center}
{\bf Abstract}
\end{center}
We give four applications of Zalcman's lemma to the dynamics of rational maps on the Riemann sphere: 
a parameter analogue of a proof of the density of repelling cycles in the Julia sets;
similarity between the Mandelbrot set and the Julia sets;
a construction of the Lyubich-Minsky lamination and its variant; 
and a unified characterization of conical points by Lyubich-Minsky 
and those by Martin-Mayer.
\end{abstract}

\section{Lemme de Zalcman}
Soient $D$ un domaine dans $\C$ et $\cF$ 
une famille d'applications holomorphes de $D$ dans $\Chat$. 
\textit{Le lemme de Zalcman} est une caract\'erisation de la (non) normalit\'e:

\begin{lem}[Lemme de Zalcman {\cite{Za}, \cite{Za2}}]\label{lem:Zalcman}
La famille $\cF$ n'est pas normale au voisinage de $z_0 \in D$ si et seulement s'il existe des suites 
$\braces{F_k}_{k \in \N} \subset \cF$,  
$\braces{\rho_k}_{k \in \N} \subset \Cstar$ avec $\rho_k \to 0$, 
et $\braces{z_k}_{k \in \N} \subset D$ avec $ z_k \to z_0$
telles que 
la suite $\psi_k(w) = F_{k}(z_k + \rho_k w)$ 
converge vers une fonction m\'eromorphe non constante
 $\psi:\C \to \Chat$ uniform\'ement sur tout compact de $\C$.
\end{lem}
On peut prendre la suite $\rho_k$ r\'eelle positive, 
mais on utilise cette version complexe dans cet article.
De plus, notons qu'on peut remplacer la suite $\psi_k$ par 
$\tilde{\psi}_k(w) = F_{k}(z_k + \rho_k w +e_k(w))$ 
o{\`u} $e_k(w) = o(\rho_k)$ sur tout compact de $\C$,
sans changer sa limite. 

\begin{figure}[htbp]
\begin{center}
\includegraphics[width=.85\textwidth]{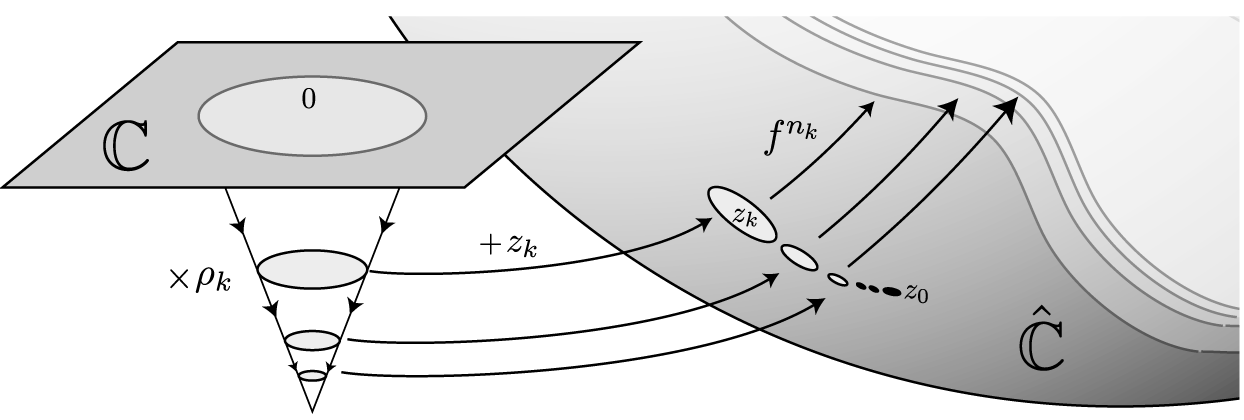}
\end{center}
\caption{Le lemme de Zalcman pour la famille $\cF = \braces{f^n}$ des it\'er\'ees d'une application rationnelle $f$.}\label{fig_zalc}
\end{figure}

\parag{Applications du lemme \`a la dynamique complexe.}
On prend une fraction rationnelle $f:\Chat \to \Chat$ 
et on applique le lemme \`a la dynamique complexe 
engendr\'ee par la famille des it\'er\'ees $\cF = \braces{f^n}_{n \ge 0}$ (voir la \figref{fig_zalc}). 
Alors le point $z_0$ dans le lemme est un point 
dans l'ensemble de Julia $J = J(f)$ de $f$.
(Dans ce cas on prend la fonction $F_k$ de la forme $F_k = f^{n_k}$ avec $n_k \to \infty$.)

Le lemme de Zalcman semble parfait pour la th\'eorie de la dynamique complexe, 
mais il y a peu de r\'esultat obtenu
en utilisant ce lemme: 
la premi\`ere application est sans doute 
une d\'emonstration simple de la densit\'e des cycles r\'epulsifs 
dans l'ensemble de Julia par Schwick \cite{Sch}.
(Elle a \'et\'e am\'elior\'ee 
par Bargmann \cite{Ba} et Berteloot-Duval \cite{BD}.
Vois aussi \cite{Za2} et \cite{Sta}.)
Steinmetz \cite{Ste} a \'etudi\'e 
des propri\'et\'es 
des fonctions m\'eromorphes engendr\'ees 
par le lemme,
en utilisant la th\'eorie de la distribution des valeurs.
Ha{\"{\i}}ssinsky \cite{Ha} et Martin-Mayer \cite{MM}
ont trouv\'e des applications 
aux ph\'enom\`enes de rigidit\'e. 
(Le lemme est aussi utilis\'e implicitement dans \cite{Mc1}.
Le principe du lemme est souvent utilis\'e dans \cite{BM}.)

Dans cet article on donnera quatre nouvelles applications du lemme de Zalcman: 
\begin{enumerate}
\item
un analogue param\'etrique de la d\'emonstration de la densit\'e 
des cycles r\'epulsifs par Schwick;
\item
une preuve alternative et simplifi\'ee 
des th\'eor\`emes de Tan Lei \cite{TL} et 
de Rivera-Letelier \cite{RL}, 
sur la ressemblance entre l'ensemble de Mandelbrot 
et les ensembles de Julia
aux param\`etres semi-hyperboliques; 
\item
une construction alternative de la lamination de Lyubich-Minsky
des fractions rationnelles  
au moyen de fonctions g\'en\'er\'ees 
par le lemme de Zalcman; et
\item
une caract\'erisation unifi\'ee
des notions de points coniques de Lyubich-Minsky et de Martin-Mayer.
\end{enumerate}

\parag{Notations.}
Dans toute la suite, $\N$ d\'esigne
l'ensemble des entiers positifs ou nul, {\it i.e.} $\N := \{0,1,2,\cdots\}$. 
Pour $x \in \C$ et $r >0$, on note $\D(x,r)$ le disque de centre $x$ et de rayon $r$. 
En particulier, $\D(0,r)$ et $\D(0,1)$ sont respectivement not\'es 
$\D(r)$ et $\D$.

Pour deux variables complexes $a$ et $b$, on \'ecrit $a \asymp b$ 
s'il existe une constante $C>1$ telle que $|a|/C \le |b| \le |a|C$.

\parag{Remerciements.}
Je voudrais remercier Jean-Yves Briend, Carlos Cabrera, 
Peter Ha{\"{\i}}ssinsky, et les rapporteurs 
dont les remarques ont permi d'am\'eliorer cet article.
Je voudrais aussi remercier le LATP pour son hospitalit\'e, 
o\`u j'ai pr\'epar\'e ce travail. 
Je voudrais adresser en particulier
 mes remerciements \`a Peter Ha{\"{\i}}ssinsky.

Ce travail est partiellement financ\'e par 
la Fondation Sumitomo, 
la Fondation Chubei Itoh, 
et la JSPS.

\section{Un analogue param\'etrique de la m\'ethode de Schwick}
On commence par un \'echauffement pour s'habituer au lemme de Zalcman.
On montre ici qu'on peut imiter l'id\'ee de Schwick \cite{Sch} 
dans l'espace param\'etrique,
et comment on peut remplacer
la m\'ethode traditionnelle utilisant le th\'eor\`eme de Montel
 par des m\'ethodes utilisant le lemme de Zalcman.
(L'origine de la m\'ethode traditionnelle 
dans l'espace param\'etrique
se trouve dans l'article \cite{Le} de Levin.)

Remarquons qu'il y a une d\'emonstration tr\`es simple 
du th\'eor\`eme de Montel en utilisant le lemme de Zalcman. 
Voir \cite{Za2} ou \cite{BM}.

\subsection{Lieu de bifurcation et param\`etres de Misiurewicz}
On consid\`ere une famille de fonctions rationnelles
param\'etr\'ees par le disque unit\'e $\D$ du plan complexe 
comme McMullen \cite{Mc2}:
soit $f:\D \times \Chat \to \Chat$
une application holomorphe de la forme $f:(t,z) \mapsto f_t(z)$,
o{\`u} $f_t:\Chat \to \Chat$ est une fraction rationnelle 
de degr\'e fix\'e $d \ge 2$.

\parag{Lieux de bifurcation et d'activit\'e.}
Le \textit{lieu de bifurcation} $B(f)$ est l'ensemble 
des param\`etres $t_0 \in \D$ tels que $t_0$ ne soit pas $J$-stable, 
{\it i.e.}, il n'y a pas de famille continue de conjugaisons
$\phi_t:J(f_t) \to J(f_{t_0})$ entre $f_t|J(f_t)$ 
et $f_{t_0}|J(f_{t_0})$ sur tout voisinage de $t_0$. 

Supposons qu'il existe une application holomorphe
$c:\D \to \Chat$ telle que $c_t:= c(t)$ soit un point 
critique de $f_t$ pour tout $t \in \D$, 
et consid\'erons la paire $(f,c)$. 
Chaque point $c_t$ est appel\'e 
le \textit{point critique marqu\'e} de $f_t$. 
On dit que $c_t$ est \textit{actif} en $t = t_0$ 
si la famille $\braces{t \mapsto f_t^n(c_t)}_{n \in \N}$
n'est pas normale sur tout voisinage de $t_0$.
Le \textit{lieu d'activit\'e} $A(f,c) \subset \D$ 
de la paire $(f,c)$ est l'ensemble des param\`etres 
$t_0 \in \D$ tels que $c_{t_0}$ est actif. 
En fait, $A(f,c) \subset B(f)$. 
Voir \cite[\S4.1]{Mc1}.

\parag{Exemple.}
Un exemple typique est la famille des polyn{\^o}mes quadratiques 
$f:(t,z) \mapsto f_t(z) = z^2 + 3t$ avec le point critique marqu\'e 
$c(t)  = 0$ pour tout $t \in \D$. On a alors $A(f,c) = B(f) = \braces{t \in \D \st 3t \in \partial M}$, o\`u $M $ est l'ensemble de Mandelbrot (voir la section 2).

\parag{Points critiques pr\'ep\'eriodiques.}
Pour une telle paire $(f,c)$, 
on note
 $\mathrm{Prep}(k,l)~(k \ge 0, l \ge 1)$ 
l'ensemble des param\`etres $t$ de $\D$ 
tels que $f_t^{k + l}(c_t) = f_t^k(c_t)$,
o{\`u} $k$ et $l$ sont minimaux. 
L'entier $k$ s'appelle le \textit{temps d'arriv\'ee} du point critique marqu\'e $c_t$.
Puisque $\mathrm{Prep}(k,l)$ est d\'etermin\'e 
par une \'equation analytique, 
il est discret ou tout le disque $\D$.
Si $k = 0$ et $t \in \mathrm{Prep}(0,l)$, 
le point critique marqu\'e $c_t$ est un point p\'eriodique 
superattractif. 
On a donc $t \notin A(f,c)$ et 
on appelle $t \in \mathrm{Prep}(0,l)$
un param\`etre {\it superattractif}.

Pour $k \ge 1$ et $l \ge 1$, 
on dit que $t \in \mathrm{Prep}(k,l)$ est 
un param\`etre de \textit{Misiurewicz} de la paire $(f,c)$ 
si $f_t^k(c_t)$ est un point p\'eriodique r\'epulsif.
On note $\mathrm{Mi}(k,l)$ 
l'ensemble des param\`etres de Misiurewicz 
dans $\mathrm{Prep}(k,l)$. 
Il est facile de montrer que 
$\mathrm{Mi}(k,l)$ 
est un sous-ensemble dense de $A(f,c)$,
par une m\'ethode simple et traditionnelle.
(Voir par exemple \cite[Thm. 2]{Le} ou \cite[Prop 2.1]{Mc2}.)

\parag{La m\'ethode de Schwick dans l'espace param\'etrique.}
On donne un peu plus d'informations sur la distribution de 
$\mathrm{Mi}(k,l)$ et $\mathrm{Prep}(k,l)$ 
lorsque $k \ge 0$ ou $l \ge 1$ est fix\'e,
en utilisant la m\'ethode de Schwick \cite{Sch}:

\begin{thm}[Distribution de Prep$(k,l)$]\label{thm_1}
Supposons que $A(f,c)$ ne soit pas vide,
et $t_0 \in \D$ soit un \'el\'ement de $A(f,c)$. 
Alors on a les propri\'et\'es suivantes:
\begin{enumerate}[\rm (1)]
\item
Si $f_{t_0}$ a un cycle r\'epulsif de p\'eriode $l \ge 3$, 
alors il existe une suite $t_j \in \mathrm{Mi}(k_j,l)$
telle que $t_j \to t_0$ et $k_j \to \infty$ quand $j \to \infty$.
\item
S'il existe un entier $k \ge 1$ avec 
$\deg (f_{t_0}, f_{t_0}^{k-1}(c_{t_0})) < \deg f_{t_0}$, 
alors 
il existe une suite $t_j \in \mathrm{Prep}(k,l_j)$
telle que $t_j \to t_0$ et $l_j \to \infty$ quand $j \to \infty$.
\item 
Pour tout $t_0 \in A(f,c)$ 
il existe une suite $t_j \in \mathrm{Prep}(0,l_j)$ 
de param\`etres superattractifs 
telle que $t_j \to t_0$ et $l_j \to \infty$ quand $j \to \infty$.
\end{enumerate}
\end{thm}
Dans (2) on note $\deg (f_{t_0}, a)$ 
le degr\'e local de $f_{t_0}$ en un point $a \in \Chat$.

Remarquons que (1) est une propri\'et\'e locale de 
$t_0 \in A(f,c)$.
Mais il garantit la densit\'e des param\`etres de Misiurewicz 
dans tout $A(f,c)$.

\begin{cor}\label{cor_B}
Supposons que le lieu de bifurcation $B(f)$ ne soit pas vide.
Alors l'ensemble des param\`etres $t$ tels que $f_t$ ait 
un point critique pr\'e-r\'epulsif est dense dans $B(f)$.
De plus, $B(f)$ est contenu dans l'adh\'erence de l'ensemble des param\`etres $t$ ayant un point critique p\'eriodique.
\end{cor}

\begin{pf}
Supposons que $t_0 \in B(f)$.
Alors $f_{t_0}$ a au moins un point critique actif
(voir \cite[Thm.4.2]{Mc1})
et donc on peut appliquer le \thmref{thm_1},
mais il peut y avoir des points critiques multiples.
Dans ce cas, on utilise  des s\'eries de Puiseaux 
(voir \cite[Prop.2.4]{Mc2} et \cite[\S 8.13]{Fo}).

Les points critiques de $\{f_{t}\}$ 
forment une vari\'et\'e analytique dans $\D \times \Chat$.
En prenant une carte locale $t \mapsto s$ 
de la forme   
$t = t_0 + \eta s^p$ pour un $|\eta| \ll 1$, 
la famille holomorphe 
$\braces{f_t \st t = t_0 + \eta s^p,~ s \in \D}$ 
repr\'esente une perturbation locale de $f_{t_0}$
telle que tout point critique 
soit param\'etr\'e holomorphiquement.
Ainsi on peut trouver un point critique marqu\'e qui est actif en $t = t_0$, et appliquer le \thmref{thm_1}.
\QED
\end{pf}

\subsection{D\'emonstration du \thmref{thm_1}}

\parag{D\'emonstration de (1).}
D'abord nous montrons (1) pour comparer 
une d\'emonstration traditionnelle
utilisant le th\'eor\`eme de Montel
avec une autre d\'emonstration
utilisant l'id\'ee de Schwick. 

Les deux d\'emonstrations suivantes 
reposent sur la d\'ependance holomorphe 
des points p\'eriodiques r\'epulsifs :
comme par hypoth\`ese il y a un cycle r\'epulsif contenant 
au moins trois points,
on peut trouver 
des applications holomorphes
$t \mapsto \braces{\al_t,\beta_t, \gamma_t}$ 
dans un petit voisinage de $t = t_0$ 
param\'etrant  trois points p\'eriodiques r\'epulsifs de $f_t$.

\parag{D\'emonstration traditionnelle.}
Comme la famille $\{F_n:t \mapsto f_{t}^n(c_t)\}_{n \in \N}$ 
n'est pas normale sur tout voisinage de $t = t_0$,
en vertu du th\'eor\`eme de Montel,
il existe des suites $t_j \to t_0$ et $n_j \to \infty$ 
telles que $f_{t_j}^{n_j}(c_{t_j}) \in \braces{\al_{t_j},\beta_{t_j}, \gamma_{t_j}}$. 
On a donc $t_j \in \mathrm{Mi}(k_j,l)$ pour un $k_j \ge 1$.

Montrons que le temps d'arriv\'ee $k_j$
n'est pas born\'e. 
Si c'\'etait le cas, il y aurait un $k \ge 1$ 
avec $f_{t}^{k + l}(c_t) = f_{t}^{k}(c_t)$ 
pour une suite $t = t_j \to t_0$,
alors $\{F_n:t \mapsto f_{t}^n(c_t)\}_{n \in \N}$ 
serait normale au voisinage de $t = t_0$. 
Mais c'est contradictoire avec l'activit\'e.  
\QED

\parag{D\'emonstration utilisant le lemme de Zalcman.}
Comme $\{F_n:t \mapsto f^n_t(c_t)\}_{n \in \N}$ n'est
pas normale en $t = t_0$, 
par le lemme de Zalcman (\lemref{lem:Zalcman}), 
il existe des suites 
$n_j \to \infty$,
$u_j \to t_0$, 
$\rho_j \to 0$ ($j \to \infty$) telles que
$\psi_j(w):=F_{n_j}(u_j+\rho_j w) 
= f_{u_j+\rho_j w}^{n_j}(c_{u_j+\rho_j w})$ 
converge vers une fonction 
m\'eromorphe non constante $\psi:\C \to \Chat$
uniform\'ement sur tout compact de $\C$.

Par le th\'eor\`eme de Picard,
il existe un $w_0 \in \C$ tel que $\psi(w_0)$ est 
un des points r\'epulsifs 
$\{\al_{t_0}, \beta_{t_0}, \gam_{t_0}\}$. 
Supposons que $\psi(w_0)=\al_{t_0}$ par exemple.
La fonction 
$w \mapsto 
f_{u_j+\rho_j w}^{n_j}(c_{u_j+\rho_j w})
-\al_{u_j+\rho_j w}$
converge vers $w \mapsto \psi(w)-\al_{t_0}$
uniform\'ement sur un disque contenant $w = w_0$, 
et donc le th\'eor\`eme de Hurwicz
affirme qu'il existe une suite $w_j \to w_0$ 
tel que 
$f_{u_j+\rho_j w_j}^{n_j}(c_{u_j+\rho_j w_j})
=\al_{u_j+\rho_j w_j}$.
Posons $t_j:=u_j+\rho_j w_j$.
On a alors $t_j \to t_0~(j \to \infty)$,
et comme $\al_{u_j+\rho_j w_j}$ est un point p\'eriodique r\'epulsif,
$t_j \in \mathrm{Mi}(k_j,l)$ pour un $k_j \ge 1$.
On peut montrer $k_j \to \infty$ comme ci-dessus.
\QED

\parag{D\'emonstrations de (3).}
Ensuite nous montrons (3) et comparons 
les deux m\'ethodes encore une fois.

\parag{Esquisse de d\'emonstration traditionnelle.}
Si $c_{t_0}$ est un point exceptionnel
({\it i.e.}, l'ensemble 
$\bigcup_{n \ge 0} f_{t_0}^{-n}(c_{t_0})$ 
est au plus deux points : voir \cite[Lem. 4.9]{Mi}), 
alors  ce point est p\'eriodique et superattractif, donc il n'est pas actif.

{\`A} l'aide d'une s\'erie de Puiseaux,
on peut trouver une autre coordonn\'ee 
au voisinage de $t_0$
et des fonctions holomorphes 
$t \mapsto \braces{a_{t}, b_{t}, c_{t}}$ 
telles que 
les graphes de ces fonctions soient disjoints
et satisfont $f_t^m(a_t) = f_t^n(b_t) =c_t$ 
pour deux entiers $n, m \ge 1$. 

Maintenant on peut appliquer
le m\^eme argument que celui de (1). 
\QED

\parag{D\'emonstration utilisant le lemme de Zalcman.}
On applique le lemme \`a la famille 
$\braces{F_n: t \mapsto f_t^n(c_t)}_{n \in \N}$. 
Supposons que $t_0 \in A(f,c)$. 
Alors on peut trouver des suites $n_j, \rho_j$ et $u_j \to t_0$ $(j \to \infty)$
telles que $\psi_j(w): = F_{n_j}(u_j + \rho_j w)$ converge vers 
une fonction m\'eromorphe non constante $\psi(w)$ 
de $\C$ uniform\'ement sur tout compact.

Comme $c_{t_0}$ n'est pas un point exceptionnel,
on peut trouver un entier $m \ge 1$ 
tel qu'il existe au moins trois solutions 
$\{ z_0, z_1, z_2\}$ de l'\'equation $f_{t_0}^m(z) = c_{t_0}$.
Supposons que 
$z_0$ n'est pas une valeur exceptionnelle 
au sens de Picard pour une telle $\psi$. 
Il existe donc un point $w_0 \in \C$ tel que 
$\psi(w_0)=z_0$ et donc $f_{t_0}^m \cc \psi(w_0) = c_{t_0}$.

Comme $\psi_j \to \psi~(j \to \infty)$ et 
$\psi_j(w) =f_{u_j + \rho_j w}^{n_j}(c_{u_j + \rho_j w})$, l'\'equation 
$$
f_{u_j + \rho_j w}^{m} \cc f_{u_j + \rho_j w}^{n_j}
(c_{u_j + \rho_j w})  =  c_{u_j + \rho_j w} 
$$
a une solution $w_j$ tendant vers $w_0$. 
Soit $t_j:= u_j + \rho_j w_j$, alors $t_j \to t_0~(j \to \infty)$ et $c_{t_j}$ est p\'eriodique de p\'eriode $l_j: = n_j+m$. 
(La p\'eriode de $c_{t_j}$ est non born\'ee, 
sinon la famille $\braces{c_t \st t \in \D}$ serait normale.) 
\QED 

\parag{D\'emonstration de (2).}
Il y a une preuve utilisant une id\'ee de Douady-Hubbard
\cite[Chapter V.2]{DH} (voir aussi \cite[Thm. 3.1]{Mc2}),
mais on doit contr{\^o}ler les points de Misiurewicz 
d\'elicatement, et la preuve n'est pas tr\`es simple.
Nous donnons ici une autre preuve utilisant le lemme de Zalcman.

Fixons le temps d'arriv\'ee $k \ge 1$. 
Par hypoth\`eses, il existe un point $a_0 \in f_{{t_0}}^{-1}(f_{{t_0}}^k(c_{{t_0}}))$ tel que 
$a_0 \neq {f_{{t_0}}^{k-1}(c_{{t_0}})}$. 
Il y a donc un entier $m \ge 1$ tel qu'il existe au moins trois solutions diff\'erentes 
$\{ z_0, z_1, z_2\}$ de l'\'equation $f_{t_0}^m(z) = a_0$. 
Prenons des suite 
$n_j \to \infty$,
$\rho_j \to 0$, et
$u_j \to t_0$ ($j \to \infty$) 
telles que 
$f^{n_j}_{u_j + \rho_j w}(u_j + \rho_j w) \to \psi(w)~(j \to \infty)$.Alors on peut trouver un $w_0 \in \C$ 
tel que $f_{t_0}^{m} \cc \psi(w_0) = a_0$ et 
$w_0$ est donc une solution de l'\'equation 
$f_{t_0}(f_{t_0}^{m} \cc \psi(w)) = f_{t_0}^k(c_{t_0})$.
Pour $j \gg 0$, il y a aussi une solution $w_j$ de l'\'equation
$$
f_{u_j + \rho_j w}( 
f_{u_j + \rho_j w}^{m} 
\cc f_{u_j + \rho_j w}^{n_j}(c_{u_j + \rho_j w})) 
 =  f_{u_j + \rho_j w}^{k}(c_{u_j + \rho_j w}) 
$$
telle que $w_j \to w_0$.
Soit $t_j := u_j + \rho_j w_j$.
Comme $t_j \to t_0$ et $f_{t_j}^{m} \cc f_{t_j}^{n_j}(c_{t_j}) \to a_0$, 
on a $f_{t_j}^{m+n_j}(c_{t_j}) \neq f_{t_j}^{k-1}(c_{t_j})$ 
et le param\`etre $t_j$ est un \'el\'ement de $\mathrm{Prep}(k, m + n_j + 1-k)$. 
Posons $l_j:= m+n_j + 1-k$.
\QED

\paragraph{Remarque.}
Schwick a aussi montr\'e la densit\'e des cycles r\'epulsifs 
dans l'ensemble de Julia pour les fonctions enti\`eres 
non polynomiales \cite{Sch},
et la m\^eme m\'ethode marche pour 
les familles de fonctions m\'eromorphes. 
Remarquons que l'analogue param\'etrique marche aussi.

\section{Ressemblance entre $M$ et $J$}
Dans cette section,
nous donnons une application du lemme de Zalcman \`a la  
dynamique des polyn{\^o}mes quadratiques.
Dans \cite{TL}, Tan Lei a montr\'e qu'au 
voisinage de tout param\`etre de Misiurewicz, 
l'ensemble de Mandelbrot 
ressemble \`a l'ensemble de Julia de ce param\`etre. 
Ce r\'esultat est am\'elior\'e par Rivera-Letelier \cite{RL}
pour les param\`etres semi-hyperboliques,
mais sa m\'ethode est diff\'erente de celle de Tan Lei.

On donne ici une d\'emonstration unifi\'ee et simplifi\'ee de ces r\'esultats, 
en utilisant le lemme du Zalcman. 
Plus pr\'ecis\'ement, nous n'utilisons pas le lemme directement,
mais nous utilisons son principe et l'id\'ee de Schwick.
De plus, on montre que 
l'ensemble de Julia du param\`etre faiblement hyperbolique
est auto-similaire au voisinage de sa valeur critique (c'est-\`a-dire de ce param\`etre).

\begin{figure}[htbp]
\begin{center}
\includegraphics[width=.75\textwidth]{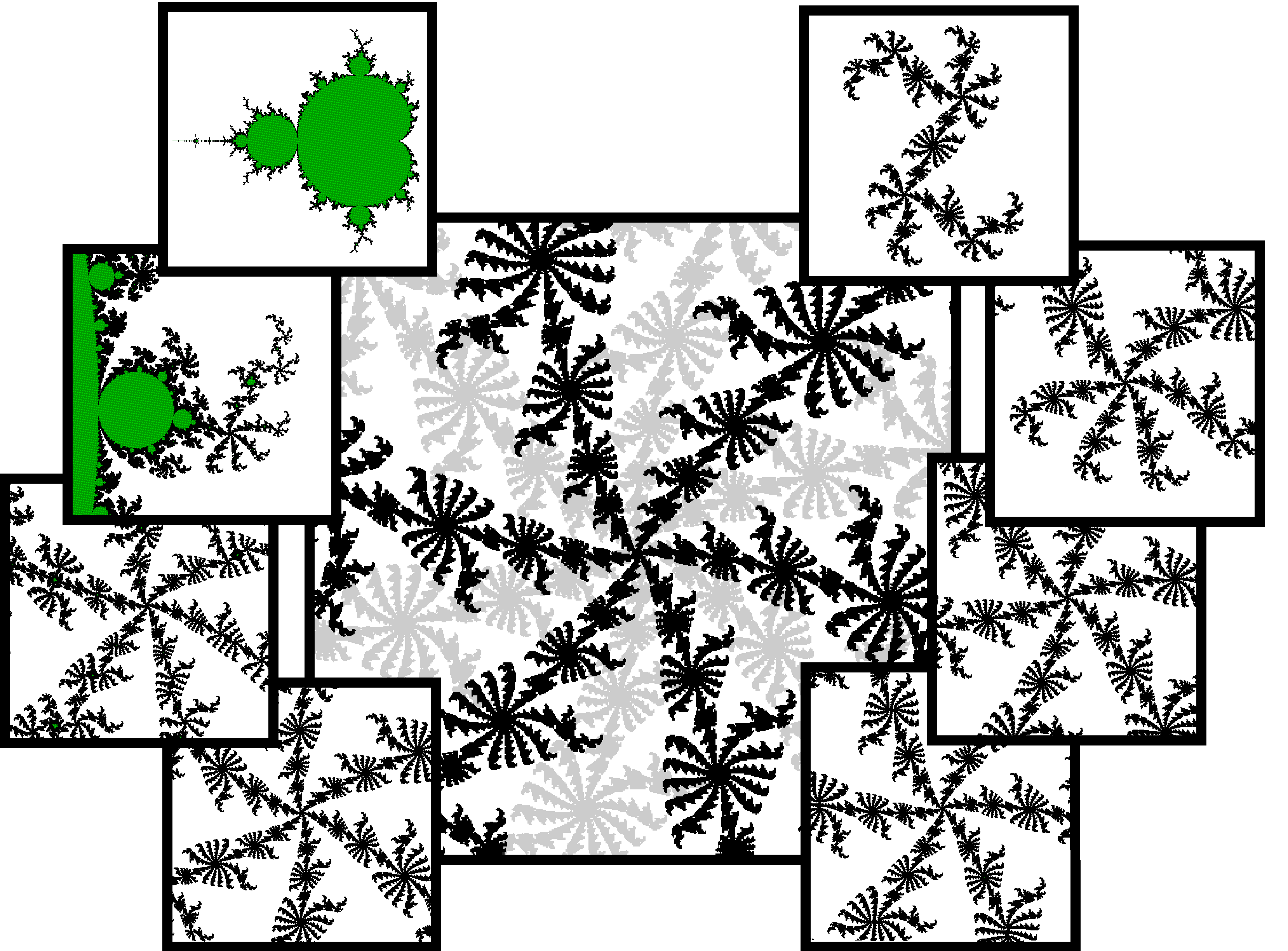}
\caption{Le tableau central repr\'esente des petites pi\`eces de 
$M$ (en gris) et de $J(f_{c_0})$ (en noir) centr\'ees en un param\`etre de Misiurewicz $c_0$ 
dans la m\^eme coordonn\'ee. 
En zoomant en arri\`ere, on voit qu'en fait ils sont  globalement diff\'erents.
}
\label{fig_MJ}
\end{center}
\end{figure}

\subsection{La ressemblance}

\parag{L'ensemble de Mandelbrot.}
Consid\'erons la famille quadratique
$$
\braces{f_c(z) = z^2 + c \st c \in \C}
$$ 
comme \cite{TL}.
Soit $M$ l'ensemble de Mandelbrot, {\it i.e.},
$$
M  :=  \braces{c \in \C : |f_c^n(c)| \le 2 ~(\forall n \in \N) }. 
$$
D'apr\`es la section 1, la fronti\`ere $\partial M$ 
est caract\'eris\'ee comme le lieu 
d'activit\'e/bifurcation du point critique $z = 0$.
Lorsqu'on a $c \in M$, on d\'efinit {\it l'ensemble de Julia rempli} de $f_c$ par 
$$
K(f_c)  :=  
\braces{z\in \C : |f_c^n(z)| \le 2 ~(\forall n \in \N) }. 
$$
En fait, on a toujours $J(f_c) = \partial K(f_c)$.

\parag{Param\`etres semi-hyperboliques.}
Un param\`etre $c_0 \in \partial M$ est appel\'e \textit{semi-hyperbolique} si l'adh\'erence de l'orbite de $c_0$ ne contient pas $c_0$ lui-m\^eme. 
Par exemple, si $c_0 \in \partial M$ est Misiurewicz, 
{\it i.e.}, $c_0$ est strictement pr\'ep\'eriodique, 
alors $c_0$ est semi-hyperbolique.

Si $c_0 \in \partial M$ est semi-hyperbolique, 
alors tous les cycles p\'eriodiques dans $\C$ sont r\'epulsifs,
et on a donc $J(f_{c_0}) = K(f_{c_0})$.
La propri\'et\'e la plus importante pour nous est le suivant:
\begin{lem}\label{lem_key}
Soit $c_0 \in \del M$ un param\`etre semi-hyperbolique. 
Alors:
\begin{enumerate}[\rm (1)]
\item
il existe des suites 
$n_k \in \N$ avec $n_k \to \infty$ et 
$\rho_k \in \Cstar$ avec $\rho_k \to 0~(k \to \infty)$ 
telles que
la suite des fonctions 
$$
\phi_k(w)  =  f_{c_0}^{n_k}(c_0 + \rho_kw)
$$ 
converge uniform\'ement sur tout compact de $\C$ 
vers une fonction $\phi:\C \to \C$ non constante; de plus,
\item
il existe une constante $Q \neq 0$ 
telle que la suite 
$$
\Phi_k(w)  :=  
f_{c_0 + Q \rho_k w}^{n_k}(c_0 + Q \rho_k w)
$$
converge vers la m\^eme fonction $\phi(w)$ 
uniform\'ement sur tout compact de $\C$.
\end{enumerate}
\end{lem}
La d\'emonstration du lemme sera donn\'ee 
dans la sous-section suivante.
Remarquons que la fonction 
$\phi(w) = \lim_{k \to \infty}  f_{c_0}^{n_k}(c_0 + \rho_kw)$ 
a la m\^eme forme que la fonction du lemme de Zalcman.

\parag{Exemple (La fonction de Poincar\'e).}
Dans le cas de Misiurewicz, 
il existe deux entiers $l, p \ge 1$
tels que $f_{c_0}^l(c_0) = f_{c_0}^{l + p}(c_0)$,
et $a_0: = f_{c_0}^l(c_0)$ est un point p\'eriodique r\'epulsif
de p\'eriode $p$.
On pose $A_0 := (f_{c_0}^l)'(c_0)$ et $\lam_0 := (f_{c_0}^p)'(a_0)$.
(On a $A_0 \neq 0$ parce que $c_0$ n'est pas un point p\'eriodique superattractif.)
Par un th\'eor\`eme de K\oe nigs 
(voir \cite[Thm.8.2, Cor.8.12]{Mi}), 
la fonction $f_{c_0}^{kp}(a_0 + w/\lam_0^{k})$
converge vers une fonction enti\`ere
$\phi:\C \to \C$ avec $\phi(\lam_0 w) = f_{c_0}^{p} \cc \phi(w)$
quand $k \to \infty$. 
La fonction $\phi$ s'appelle {\it la fonction de Poincar\'e}.
 Comme on a $f_{c_0}^l(c_0 + \Delta z) = a_0 + A_0 \Delta z +o(\Delta z)$ pour tout $\Delta z \approx 0$, 
$$
\phi(w) = \lim_{k\to \infty} 
f_{c_0}^{kp}\paren{a_0 + \frac{w}{\lam_0^{k}}}
 =  
\lim_{k\to \infty} f_{c_0}^{l + kp}\paren{c_0 + \frac{w}{A_0\lam_0^{k}} + o(\lam_0^{-k})}
$$
sur tout compact de $\C$.
On obtient donc $\phi(w) = \lim_kf_{c_0}^{l + kp}(c_0 + w/(A_0\lam_0^{k}))$, et alors on peut prendre $n_k = l + kp$ et 
$\rho_k = 1/(A_0\lam_0^{k})$ dans le \lemref{lem_key}.

\parag{Topologie de Hausdorff.}
On rappelle la \textit{topologie de Hausdorff} 
dans l'espace $\mathrm{Comp}^\ast(\C)$ des sous-ensembles compacts (non vides) de $\C$. 

Pour une suite $\{K_k\}_{k\in \N} \subset \mathrm{Comp}^\ast(\C)$, 
on dit que $K_k$ converge vers $K \in \mathrm{Comp}^\ast(\C)$ 
quand $k \to \infty$ si pour tout $\e>0$, 
il existe $k_0 \in \N$ tel que 
$K \subset \mathrm{N}_\e(K_k)$ et $K_k \subset \mathrm{N}_\e(K)$ 
pour tout $k \ge k_0$, 
o{\`u} $\mathrm{N}_\e(\cdot)$ est le $\e$-voisinage ouvert dans $\C$.

Posons $\D(r): = \braces{|z| < r}$.
Pour $K \subset \C$ ferm\'e, 
notons $[K]_r$ l'ensemble $(K \cap \D(r))\cup \partial \D(r) 
\in \mathrm{Comp}^\ast(\C)$.
Pour des constantes $a \in \Cstar$ et $b \in \C$, 
notons $a(K-b)$ l'ensemble $\braces{a(z-b) \st z \in K}$.

\parag{Ressemblance.}
Soit $c_0$ un param\`etre semi-hyperbolique.
Par le \lemref{lem_key}, on peut trouver 
des suites $n_k \to \infty$ et $\rho_k \to 0$
telles que $\phi_k(w) = f_{c_0}^{n_k}(c_0 + \rho_k w)$ 
converge vers une fonction enti\`ere non constante, 
$\phi:\C \to \C$.
(Remarquons que $\phi$ n'a pas de p{\^o}les, 
comme $\phi_k$ est un polyn\^ome.)
On a aussi une constante $Q \neq 0$ telle que
$\Phi_k(w) = f_{c_0 + Q\rho_k w}^{n_k}(c_0 + Q\rho_k w)$
converge vers la m\^eme fonction $\phi:\C \to \C$.

Le r\'esultat principal de cette section est la modification suivante de  \cite{TL, RL}:

\begin{thm}[Ressemblance entre $M$ et $J$]\label{thm_TL}
Soient $c_0 \in \del M$ semi-hyperbolique
et $\phi:\C \to \C$ la fonction limite de 
$\phi_k(w) = f_{c_0}^{n_k}(c_0 + \rho_k w)$
et de $\Phi_k(w) = f_{c_0 + Q\rho_k w}^{n_k}(c_0 + Q\rho_k w)$ 
ci-dessus. 
Posons $\cJ: = \phi^{-1}(J) \subset \C$, 
o{\`u} $J = J(f_{c_0})$ est l'ensemble de Julia de $f_{c_0}$.
Alors pour tout $r>0$, 
\begin{enumerate}[\rm (a)]
\item 
$\gauss{\rho_k^{-1}(J-c_0)}_r  \to [\cJ]_r~~(k \to \infty)$
\item 
$\gauss{Q^{-1} \rho_k^{-1}(M-c_0)}_r  \to [\cJ]_r~~(k \to \infty)$
\end{enumerate}
dans la topologie de Hausdorff.
\end{thm}

\paragraph{Remarque.}
Rivera-Letelier a montr\'e que 
le distance de Hausdorff entre 
$[\lam (J-c_0)]_R$ et $[M-c_0]_R$
est $O(R^{3/2})$ quand $R \to 0$,  
o\`u $\lam=\sum_{n \ge 0}1/(f_{c_0}^n)'(c_0)$. 
Voir le premier corollaire de \cite[p.290, Thm. B]{RL}. 
(Dans l'in\'egalit\'e du corollaire $Cr^{1/d}$ est en fait $Cr^{1+1/d}$.)  
Donc on a $Q=\lam$ et le distance de Hausdorff entre 
$[\rho_k^{-1}(J-c_0)]_r$ et $[Q^{-1} \rho_k^{-1}(M-c_0)]_r$
est $O(r^{3/2}\sqrt{\rho_k})$ pour tout $r$ fix\'e.

\paragraph{D\'emonstration de (a).}
Comme $f_{c_0}^{n}(J) = J$, 
on a $\gauss{\rho_k^{-1}(J-c_0)}_r = \gauss{\phi_k^{-1}(J)}_r$. 
Par d\'efinition de $\cJ$, 
on a $[\cJ]_r = \gauss{\phi^{-1}(J)}_r$. 
Puisque $\phi_k$ converge vers $\phi$ uniform\'ement 
dans $\overline{\D(r)}$, 
on obtient (a).

\paragraph{D\'emonstration de (b).}
Posons $\cM_k:= Q^{-1} \rho_k^{-1}(M-c_0)$.
Nous montrons que pour tout $\e >0$, 
$$
\gauss{\cM_k}_r ~\subset~ \mathrm{N}_\e([\cJ]_r)
~~\textit{et}~~~
 [\cJ]_r~\subset~ \mathrm{N}_\e(\gauss{\cM_k}_r)
$$ 
quand $k$ est assez grand.

D'abord, l'ensemble 
$\overline{\D(r)}-\mathrm{N}_\e(\cJ)$ est compact,
donc il existe un entier $N = N(\e)$ tel que
$|f_{c_0}^N \cc \phi(w)| >2$ pour tout $w \in \overline{\D(r)}-\mathrm{N}_\e(\cJ)$.
Comme $\Phi_k(w)$ converge vers $\phi(w)$ uniform\'ement sur tout 
compact de $\C$ (\lemref{lem_key}), on a 
$$
\abs{f_{c_0 + Q \rho_k w}^{N + n_k}(c_0 + Q \rho_k w)} ~>~2
$$
pour tout $k \gg 0$.
Par cons\'equent, on a $c_0 + Q \rho_k w \notin M$, {\it i.e.}, $w \notin \cM_k$. Du coup, l'inclusion
$
\gauss{\cM_k}_r \subset \mathrm{N}_\e([\cJ]_r)
$ 
est v\'erifi\'ee.

Pour l'autre inclusion, prenons 
un ensemble fini $E \subset [\cJ]_r$
tel que le $\e/2$-voisinage de $E$ contienne $[\cJ]_r$.
Il suffit de trouver 
une suite $w_k \in [\cM_k]_r$ pour chaque $w_0 \in E$
telle que 
$|w_0-w_k| < \e/2$ pour tout $k \gg 0$.

Posons $\Delta: = \D(w_0,\e/2)$.
Dans le cas $\Delta \cap \del \D(r) \neq \emptyset$,
on prend la suite $w_k$ dans $\del \D(r)$.

Dans le cas $\Delta \subset \D(r)$,
puisque $\phi(w_0) \in J(f_{c_0})$ et 
les cycles r\'epulsifs sont denses dans $J(f_{c_0})$,
on peut trouver un $w_0'$ tel que $|w_0-w_0'|< \e/4$ 
et tel que $\phi(w_0')$ soit un point p\'eriodique r\'epulsif de p\'eriode $m$.
Alors $w = w_0'$ est un z\'ero de la fonction $\chi:w \mapsto f_{c_0}^m(\phi(w))-\phi(w)$.
Consid\'erons la fonction $\chi_k: w \mapsto f_{c_0 + Q \rho_k w}^m(\Phi_k(w))-\Phi_k(w)$.
La fonction $\Phi_k$ converge vers $\phi$ 
uniform\'ement sur tout compact de $\C$,
donc pour $k \gg 0$ il existe un z\'ero 
$w_k \in \Delta$ de $\chi_k$ avec $|w_k - w_0'| < \e /4$.
En particulier,
$c_k:= c_0 + Q \rho_k w_k$ satisfait 
$f_{c_k}^{n_k + m}(c_k) = f_{c_k}^{n_k}(c_k)$.
On a donc $c_k \in M$.
Par suite, on peut trouver une suite $w_k \in \cM_k$ avec $|w_k -w_0|<\e/2$.
\hfill \QED

\subsection{D\'emonstration du \lemref{lem_key}}

\parag{Ensemble hyperbolique.}
Soit $c_0 \in \partial M$ un param\`etre semi-hyperbolique, et 
$X_0$ l'ensemble des points d'adh\'erence de l'orbite de $c_0$. 
L'ensemble $X_0$ est un \textit{ensemble hyperbolique}, 
c'est-\`a-dire, $X_0$ est compact, $f_{c_0}(X_0) = X_0$, 
et il existe des constantes $\kappa, \eta >0$ telles que 
$|(f_{c_0}^n)'(x)| \ge \kappa (1 + \eta)^n$ 
pour tout $x \in X_0$ et tout $n \in \N$ (voir \cite{CJY}).

Par exemple, si $c_0 \in \del M$ est Misiurewicz, 
alors l'orbite $\braces{f_{c_0}^n(c_0)}_{n \in \N}$ tombe en un temps fini sur un cycle r\'epulsif. Dans ce cas, $X_0$ est ce cycle r\'epulsif.
En fait, l'orbite de $c_{0}$ toujours tombe sur $X_{0}$:

\begin{lem}\label{lem_X_0}
Si $c_0 \in \partial M$ est semi-hyperbolique,
il existe un $l \in \N$ tel que $f_{c_0}^l(c_0) \in X_0$. 
\end{lem}

\parag{D\'emonstration.}
Par hyperbolicit\'e de $X_0$, on peut prendre un $p \in \N$ tel que $|(f_{c_0}^p)'(x)| \ge 3$ pour tout $x \in X_0$. 
Comme $X_0$ est compact, il existe un $\delta>0$ 
tel que si $z \in \mathrm{N}_\delta(X_{0})-X_0$ on a 
$\mathrm{dist}(f_{c_0}^{p}(z),X_0) \ge 2\, \mathrm{dist}(z,X_0)$.
Si $f_{c_0}^{n}(c_{0})\notin X_0$ pour tout  $n \in \N$,
il y a un autre point d'adh\'erence de l'orbite de $c_0$
dans $\Chat - \mathrm{N}_\delta(X_{0})$, 
mais c'est contradictoire avec la d\'efinition de $X_{0}$.
\QED{\small (\lemref{lem_X_0})}

Pour montrer le \lemref{lem_key}, on utilise une id\'ee de \cite[Lemma 4.7]{LM}:

\parag{D\'emonstration du \lemref{lem_key}, (1).}
On fixe $f:= f_{c_0}$ avec $c_0 \in \partial M$ 
semi-hyperbolique.
Par le \lemref{lem_X_0}, il existe un $l >0$ 
tel que $a_0: = f^l(c_0) \in X_0$. 
Pour $m \in \N$, posons $a_m:=f^m(a_0) \in X_0$,
$\lam_m := (f^m)'(a_0)$, et $T_m(w):= a_0 + w/\lam_m~(w \in \C)$.
Puisque $|\lam_m| \ge \kappa(1 + \eta)^m \to \infty~(m \to \infty)$,
on peut trouver un $p \in \N$ et un $\mu_0>1$ tels que 
$|(f^{p})'(x)| \ge \mu_0$ 
pour tout $x \in X_0$. 
Comme $X_0$ est l'ensemble $\omega$-limite de $0$ et de $a_0$,
il existe un $x_0 \in X_0$ et une suite $\braces{m(k)}_{k \ge 0} \subset p\N$ tels que $a_{m(k)} \to x_0$ quand $k \to \infty$. 
(Par exemple, si $c_0$ est Misiurewicz et 
$c_0$ tombe sur un point p\'eriodique r\'epulsif 
$a_0$ de p\'eriode $p \ge 1$,
alors on peut choisir $x_0 = a_0$, $\mu_0 = |(f^p)'(a_0)|$, et $m(k) = kp$.)
Soit $\psi_k:= f^{m(k)}\cc T_{m(k)}:\C \to \C$. 
Alors on a:

\begin{prop}\label{prop_Cauchy}
La suite $\braces{\psi_k}_{k\ge 0}$ converge uniform\'ement
sur tout compact de $\C$ vers une application enti\`ere $\psi$ 
telle que $\psi'(0) = 1$ et $\psi(0) = x_0$. 
\end{prop}

\parag{D\'emonstration.}
D'abord on montre qu'il existe un $\delta >0$ tel que 
$\psi_k|\D(\delta)$ est univalente pour tout $k \in \N$:
comme $|(f^p)'|\ge \mu_0>1$ sur $X_0$, il existe un $\delta_0 > 0$ tel que 
$f^p|{\D(x,\delta_0)}$ soit univalente et 
$\D(f^p(x),\delta_0) \Subset f^p(\D(x,\delta_0))$ pour tout $x \in X_0$.
Il existe donc une application inverse $h_k: \D(a_{m(k)},\delta_0) \to \C$ de $\psi_k = f^{m(k)} \cc T_{m(k)}$ 
telle que $h_k$ soit univalente et $h_k(a_{m(k)}) =h_k'(a_{m(k)})-1 = 0$. 
Par le th\'eor\`eme du quart de Koebe, 
$\D(\delta_0/4) \subset h_k(\D(a_{m(k)},\delta_0))$. 
Soit $\delta := \delta_0/4$. 
Alors la famille $\braces{\psi_k|\D(\delta)}_{k \ge 0}$ 
est univalente, donc uniform\'ement born\'ee sur tout l'ensemble compact par un th\'eor\`eme de Koebe (voir \cite[\S 2.3]{Du}). 
En particulier, la famille est \'equicontinue.

Ensuite on montre que la famille $\braces{\psi_k}_{k \ge 0}$
est normale sur $\C$: 
soit $G_{N,j}: = T_{m(N)}^{-1} \cc 
f^{m(N + j)-m(N)} \cc T_{m(N + j)}$ pour que $\psi_{N + j} = \psi_{N} \cc G_{N,j}$. 
En utilisant le th\'eor\`eme du quart de Koebe encore une fois, 
on voit que $G_{N,j}$ est univalente sur $\D(\delta |\lam_{m(N)}|)$.
Posons 
$$
w_{N,j}: = G_{N,j}(0) = h_N \cc \psi_{N + j}(0) = h_N(a_{m(N + j)}).
$$ 
Comme $h_N$ est univalente sur $\D(a_{m(N)},\delta_0)$ 
avec $h_N'(a_{m(N)}) = 1$ pour tout $N$, 
et comme $\braces{a_{m(k)}}_{k \ge 0}$ est de Cauchy,  
on obtient  $w_{N,j} \to 0$ quand $N \to \infty$.
Soit $\tilde{G}_{N,j}(w):= G_{N,j}(w)-w_{N,j}$. 
Alors $\tilde{G}_{N,j}$ est univalente sur $\D(\delta |\lam_{m(N)}|)$,
et on a $\tilde{G}_{N,j}(0) = \tilde{G}_{N,j}'(0)-1 = 0$.

On prend un $r>0$ arbitrairement grand.
Par un th\'eor\`eme de Koebe (voir \cite[\S 2.3]{Du}),
on a $|\tilde{G}_{N,j}'(w)-1| \le C |w|/|\lam_{m(N)}|$ sur $\D(r)$, 
o\`u $C = C(\delta, r)$ est une constante ind\'ependante de $N$ assez grand.
Alors on a $\tilde{G}_{N,j}(w) = w + O(r^2/|\lam_{m(N)}|)$, 
donc $G_{N,j}(w) = w + O(r^2/|\lam_{m(N)}|) + 
O(|a_{m(N + j)}-a_{m(N)}|)$.
En particulier, si $N$ est assez grand, 
on a ${G}_{N,j}(\D(r)) \subset \D(2r)$.
Soit $R = R(r,N) = \max\{|\psi_N(w)| \st |w| \le 2r\}$.
Alors 
$$\psi_{N + j}(\D(r)) = \psi_{N} (G_{N,j}(\D(r)))
\subset \D(R)
$$ 
pour tout $j \in \N$.
Puisque la famille $\{\psi_k\}_{k \ge 0}$ est 
uniform\'ement born\'ee sur tout compact de $\C$,
la famille $\{\psi_k\}_{k \ge 0}$ est normale sur $\C$.

Supposons que $w \in \D(\delta)$. 
Alors $G_{N,j}(w) \to w$ uniform\'ement sur $\D(\delta)$ 
quand $N \to \infty$.
Comme $|\psi_{N,j}(w)-\psi_{N}(w)|
 = |\psi_{N}(G_{N,j}(w))-\psi_{N}(w)| \to 0$ quand $N \to \infty$  
 et comme $\braces{\psi_N}$ est \'equicontinue 
 pour tout $N$ sur $\D(\delta)$,
la suite $\braces{\psi_{k}|\D(\delta)}_{k \ge 0}$ 
est de Cauchy. 
Ainsi la limite $\psi = \lim \psi_k$ avec $\psi'(0) = 1$ et $\psi(0) = x_0$ existe, et elle est enti\`ere par normalit\'e de $\braces{\psi_{k}}_{k \ge 0}$ sur $\C$. 
\QED{\small (\propref{prop_Cauchy})}

\parag{} 
Posons 
$n_k := m(k) + l$, $\rho_k:= 1/(f^{n_k})'(c_0)$.
$\phi_k(w):=f^{n_k}(c_0 + \rho_k w)$, 
et $A_0:= (f^l)'(c_0) \neq 0$. 
On a alors $A_0 \lam_{m(k)} = 1/\rho_k$ et 
$$
f^l(z)  =  a_0 + A_0(z-c_0) +  o(z-c_0),
$$
donc 
$$
\phi_k(w)  =  f^{m(k)}(a_0 + A_0 \rho_k w + o(\rho_k))
~\to~ \lim \psi_{k}(w)  =  \psi(w)
$$
uniform\'ement sur tout compact de $\C$.
Soit $\phi: = \psi$. 
\QED{\small(\lemref{lem_key},(1))}

\medskip

Pour montrer le \lemref{lem_key}(2), 
on utilise la stabilit\'e dynamique de l'ensemble hyperbolique $X_0$
(voir \cite[\S 1]{Shi}):

\begin{prop}[Mouvement holomorphe de $X_0$]\label{prop_holo}
Il existe un voisinage $U$ de $c_0$ dans $\C$
tel qu'on ait un mouvement holomorphe dynamique $\chi:X_0\times U \to \C$:  
\begin{itemize}
\item
pour tout point $x \in X_0$, $\chi(x, c_0) = x$;
\item
pour tout point $x \in X_0$ fix\'e, 
l'application $\chi^{x}: c \mapsto \chi(x,c)$ est holomorphe; et  
\item
pour tout param\`etre $c \in U$ fix\'e, 
l'application $\chi_c:x \mapsto \chi(x,c)$ est quasiconforme
telle que $f_c \cc \chi_c = \chi_c \cc f_{c_0}$ sur $X_0$.

\end{itemize}
\end{prop}
Par exemple, si $c_0 \in \del M$ est Misiurewicz, 
$\chi_c(X_0)$ est un cycle r\'epulsif de $f_c$.
Dans ce cas, 
la fonction de Poincar\'e de chaque point 
p\'eriodique dans $\chi_c(X_0)$ d\'epend holomorphiquement de $c$. 
Plus g\'en\'eralement, si $c_0 \in \partial M$ est 
semi-hyperbolique, on a:

\begin{lem}\label{lem_holo_dep}
Pour $c \in U$, le voisinage de $c_0$ ci-dessus, 
posons $a(c): = \chi_c(a_0)$, $\lam_m(c):= (f_c^m)'(a(c))~(m \in \N)$, et soit $\braces{m(k)}_{k \in \N} \subset p\N$ la suite 
trouv\'ee dans la d\'emonstration du \lemref{lem_key}, 
{\rm (1)} telle que
$\phi(w) = \lim_{k} f_{c_0}^{m(k)}(a_0 + w/\lam_{m(k)}(c_0))$. 
Alors la suite $\psi_{m(k)}^{c}(w):= f_c^{m(k)}(a(c) + w/\lam_{m(k)}(c))$ converge vers une fonction m\'eromorphe $\phi^c(w)$
uniform\'ement sur tout compact de $\C$.
En plus, on a $\phi^{c_0} = \phi$,
et si on fixe $w \in \C$, 
l'application $c \mapsto \phi^c(w)$ est holomorphe en $c = c_0$.
\end{lem}

\parag{D\'emonstration.}
On peut supposer que le voisinage $U$ est assez petit tel que
$\lam_{m}(c) \ge \kappa(1 + \eta)^m$ 
pour tout $c \in U$ et $m \in \N$. 
Soit $\braces{m(k)}_{k \in \N}$ la suite de 
la d\'emonstration du \lemref{lem_key}(1)
 avec $a_{m(k)} \to x_0~(k \to \infty)$.
 Par la conjugaison $f_c \cc \chi_c = \chi_c \cc f_{c_0}$ sur $X_0$ de \propref{prop_holo}, on a aussi 
 $f_c^{m(k)}(a(c)) \to \chi_c(x_0)~(k \to \infty)$ 
 sur $\chi_c(X_0)$.
De la m\^eme fa\c{c}on que dans 
la d\'emonstration du \lemref{lem_key}(1), 
on peut montrer la convergence de 
$\{\phi_{m(k)}^{c}\}_{k \in \N}$ sur $\C$
et que la fonction $c \mapsto \phi^{c}(w)$
est holomorphe sur $U$ pour chaque $w \in \C$.
\QED{\small (\lemref{lem_holo_dep})}

\parag{D\'emonstration du \lemref{lem_key}, (2).}
On pose $c:= c_0 + Q \rho_k w$ avec une constante $Q \in \Cstar$.
On pose aussi $\Phi_k(w) := f_{c}^{n_k}(c) = f_{c}^{m(k) + l}(c)$
et $b(c):= f_c^{l}(c)$. 
Par \cite[Appendix 2]{RL} (ou \cite[Thm.1.1]{vS}),
on a 
$b(c)-a(c) = B_0(c-c_0) +o(c-c_0)
$ avec $B_0 \neq 0$.
En utilisant $\rho_k^{-1} = A_0 \cdot \lam_{m(k)}(c_0)$,
 on a ainsi
$$
b(c)  =  a(c) + B_0Q \rho_k w +o(\rho_k) 
 =  a(c) + 
\frac{B_0Q}{A_0} 
 \cdot 
\frac{ \lam_{m(k)}(c) }{ \lam_{m(k)}(c_0)} 
 \cdot \frac{w}{\lam_{m(k)}(c)}
+o(\rho_k). 
$$
On prend ici $Q:= A_0/B_0$. 
Comme $X_0$ est compact avec $0 \notin X_0$, on a $|f_c'(\chi_c(x))/f_{c_0}'(x)-1| = O(c-c_0)$ uniform\'ement pour tout $x \in X_0$, 
donc quand $m = m(k) \to \infty$, on obtient
$$
\log \frac{ \lam_{m}(c) }{ \lam_{m}(c_0)} 
 =  \sum_{i = 0}^{m-1} 
 \log\frac{f_c'(f_c^i(a(c))}{f_{c_0}'(f_{c_0}^i(a({c_0}))}
 =  m \cdot O(c-c_0)
 =  O\paren{\frac{m}{\lam_{m}(c_0)}}
  ~\to~ 0.
$$
(Rappelons que $X_0$ est hyperbolique, puis $|\lam_{m}(c_0)| \ge \kappa(1+\eta)^m$.)
Puisque $\Phi_k(w) =f_c^{m(k)}(b(c))$ et que $\lim_{c \to c_0} \phi^c(w) = \phi(w)$ uniform\'ement sur tout compact de $\C$,
on obtient
$$
\lim_{k \to \infty} \Phi_k(w) 
 =  \lim_{k \to \infty} 
f_c^{m(k)}\paren{a(c) +  \frac{w}{\lam_{m(k)}(c)}+o(\rho_k)}
 =  \phi(w).
$$ 
\QED

\medbreak
\parag{Remarque.}
Si on d\'emontre les convergences $\phi_k \to \phi$ et $\Phi_k \to \phi$ juste sur un disque, la d\'emonstration est plus facile. En fait, $\phi'_k(0) = 1$ et on peut utiliser la normalit\'e des fonctions univalentes sur $\D(\delta)$. 
On a alors une version faible du \lemref{lem_key}(1).
Remarquons que la \propref{prop:MM-conical} (\cite{Ha}, \cite{MM}) 
est une g\'en\'eralisation de cette version faible.

\subsection{Auto-similitude de $J$ aux param\`etres faiblement hyperboliques}
Nous donnons ici une g\'en\'eralisation
du \thmref{thm_TL}(a).
D'apr\`es Ha\"{\i}ssinsky \cite{Ha},
on dit que $c_0 \in \partial M$ est 
\textit{faiblement hyperbolique}
s'il existe un $d \in \N$ et un $r >0$ tels que 
pour $D_n: = \D(f_{c_0}^n(c_0),r)$
et la composante connexe $D_n'$ de 
$f_{c_0}^{-n}(D_n)$ contenant $c_0$, 
$\deg(f_{c_0}^n:D_n' \to D_n) \le d$ pour un nombre infini de $n \in \N$.

\begin{thm}[Auto-similitude de $J$]\label{thm_weakly_hyp}
Soit $c_0 \in \del M$ faiblement hyperbolique.
Alors il existe une fonction holomorphe non constante $\phi:\D \to \C$ avec $\phi(0) \in J = J(f_{c_0})$ et une suite $\rho_k \in \Cstar$ avec $\rho_k \to 0$ telles que pour tout $0<r<1$, 
$$
\gauss{\rho_k^{-1}(J-c_0)}_r  ~\to~ [\phi^{-1}(J)]_r~~(k \to \infty)
$$
dans la topologie de Hausdorff. 
\end{thm}

\begin{pf}
Par la \propref{prop:MM-conical} (\cite{Ha}, \cite{MM}) de la section 4, il existe des suites $n_k \in \N$ avec $n_k \to \infty$ et $\rho_k \in \Cstar$ avec $\rho_k \to 0$ telles que 
le polyn\^ome $\phi_k(w) = f_{c_0}^{n_k}(c_0 + \rho_k w)$ converge vers
une fonction holomorphe non constante $\phi:\D \to \C$ 
uniform\'ement sur tout compact de $\D$. 
En particulier, $\phi(0) = \lim_k f_{c_0}^{n_k}(c_0) \in J$, 
donc $\phi^{-1}(J)$ n'est pas vide. 
Du coup on peut appliquer la d\'emonstration du \thmref{thm_TL}(a).

\QED
\end{pf}

\section{Lamination de Lyubich-Minsky}
Dans cette section, nous donnons une construction de la lamination de Lyubich-Minsky par le lemme de Zalcman. 

Pour une fraction rationnelle $f:\Chat \to \Chat$, 
M.~Lyubich et Y.~Minsky ont d\'efini une lamination $\cA_f$ par surfaces de Riemann.
\`A partir de $\cA_f$, ils ont construit 
une lamination par 3-vari\'et\'es hyperboliques
comme un analogue des 3-vari\'et\'es hyperboliques
des groupes kleiniens,
et ils ont d\'emontr\'e un th\'eor\`eme de rigidit\'e 
en utilisant des m\'ethodes de g\'eom\'etrie hyperbolique. 
(Voir \cite{LM} ou \cite[Chap. 3]{KL}.)

L'ingr\'edient principal de leurs laminations
est un ensemble de fonctions m\'eromorphes
g\'en\'er\'ees par la dynamique inverse de $f$.
On montre ici que si $f$ v\'erifie une condition (plus g\'en\'erale que la parabolicit\'e), 
on peut utiliser l'ensemble des fonctions m\'eromorphes 
g\'en\'er\'ees par le lemme de Zalcman
pour la famille $\cF:= \braces{f^n}_{n \in \N}$.

\subsection{Construction de la lamination de Lyubich-Minsky}
Soit $f:\Chat \to \Chat$ une fraction rationnelle de degr\'e 
au moins deux.
On commence par la construction de la lamination $\cA^{\cLM}$ ($ = \cA_f$ ci-dessus) due \`a Lyubich et Minsky.

\parag{Les orbites inverses r\'eguli\`eres.}
On dit qu'une suite $\zhat = (z_0,z_{-1},z_{-2} \cdots)$ dans 
le produit infini $\Chat^{\N}$
est une \textit{orbite inverse} de $f$
si $f(z_{-n - 1}) = z_{-n}$ pour tout $n \in \N$.
On dit qu'une orbite inverse $\zhat = (z_{-n})$ de $f$ est 
\textit{r\'eguli\`ere} s'il existe un voisinage $U_0$ de $z_0$
tel que la suite des ant\'ec\'edents it\'er\'es 
$\cdots \to U_{-2} \to U_{-1} \to U_{0}$ 
contenant l'orbite inverse 
$\cdots \mapsto z_{-2} \mapsto z_{-1} \mapsto z_{0}$ 
soit univalente sur $U_{-n}$ pour tout $n$ assez grand. 

Soit $\cR \subset \Chat^\N$ l'ensemble de toutes les orbites inverses r\'eguli\`eres. 
Une composante connexe par arcs de $\cR$ est appel\'ee 
\textit{une feuille}.  
La feuille contenant $\zhat$ est not\'ee  $L(\zhat)$. 
Les feuilles de $\cR$ sont isomorphes \`a $\C$, $\D$, ou \`a un anneau (\cite[Lemma 3.3, Lemma 3.4]{LM}).

Soit $\cR_\C$ l'ensemble de toutes les  feuilles $L(\zhat)$ de $\cR$ 
isomorphes \`a $\C$. Par exemple, si la suite $(z_{-n})$ est donn\'ee par un cycle r\'epulsif, on a $L(\zhat) \simeq \C$. 
Si $L(\zhat) \simeq \C$, il existe une \textit{uniformisation}
$\phi:\C \to L(\zhat)$ (\cite[\S 4]{LM}).

\parag{Fonctions de Lyubich-Minsky.}
Soit $\cU$ l'espace de toutes les fonctions m\'eromorphes 
non constantes de $\C$.
(La topologie de $\cU$ est celle de la convergence uniforme sur les compacts de $\C$.)
Nous consid\'erons la classe suivante de fonctions dans $\cU$
qui est tr\`es importante dans toute la suite.

Soit $\pi_0:\zhat = (z_{-n}) \mapsto z_0$ la projection sur $\Chat$.
Lorsqu'on a une uniformisation $\phi:\C \to L(\zhat)$, 
l'application $\psi = \pi_0 \cc \phi$ est une fonction 
m\'eromorphe non constante de $\C$; 
{\it i.e.,} $\psi$ est un  \'el\'ement  de $\cU$. 
Pour une orbite inverse $\zhat = (z_{-n}) \in \cR_\C$ de $f$, 
on pose 
$$
\cLM(\zhat)  :=  \braces{\psi =\pi_0 \cc \phi \in \cU 
\st \exists \phi:\C \to L(\zhat) \text{~uniformisation t.q.~} \phi(0)=\zhat}
$$
et $\cLM := \bigcup_{\zhat \in \cR_\C} \cLM(\zhat)$.
On dit qu'un \'el\'ement $\psi \in \cLM \subset \cU$ est 
une \textit{fonction de Lyubich-Minsky} 
(ou \textit{LM-fonction}).

On note l'ensemble des applications affines (complexes) de $\C$ 
par $\Aff$. 
Alors il n'est pas difficile de s'assurer que 
$\cLM$ satisfait 
$$
f \cc \cLM = \cLM =\cLM \cc \Aff.
$$
Plus pr\'ecis\'ement, on a: 
\begin{enumerate}
\item
Si $\psi \in \cLM$, alors $f \cc \psi \in \cLM$ et il existe
une $\psi_1 \in \cLM$ telle que $\psi = f \cc \psi_1$. 
\item 
Si $\psi \in \cLM$ et $\delta:\C \to \C$ est une application affine, 
alors $\psi \cc \delta \in \cLM$.
\end{enumerate}

\parag{Construction de la LM-lamination.}
Pour construire la lamination de Lyubich-Minsky, 
on consid\`ere plus g\'en\'eralement 
un sous-ensemble $\cK \neq \emptyset$ de $\cU$ avec
$$
f\cc \cK = \cK = \cK \cc \Aff.
$$
On d\'efinit alors une lamination $\cA^\cK$ comme suite:
soit $\widehat{\cU}:= \cU^{\N}$ l'ensemble des suites 
de la forme
$$
\psihat :=  (\psi_0,\psi_{-1},\psi_{-2},\ldots)
$$
o\`u $\psi_{-n} \in \cU$ pour tout $n \in \N$.
Notons $\widehat{\cK}$ l'ensemble des \'el\'ements 
$\psihat = (\psi_{-n})$ de $\widehat{\cU}$
tels que $\psi_{-n} \in \cK$ et  
$f \cc \psi_{-n} =\psi_{-n + 1}$ pour tout $n \in \N$. 
(En utilisant la propri\'et\'e $f \cc \cK = \cK$,
on a $\widehat{\cK}\neq\emptyset$.)
Alors pour tout $\psihat = (\psi_{-n}) \in \widehat{\cK}$ il existe une orbite inverse $\zhat = (z_{-n})$ telle que 
$\psi_{-n}(0) = z_{-n}$ pour tout $n \in \N$.

On d\'efinit les relations d'\'equivalence suivantes
sur $\cU$ et $\widehat{\cU}$: 
on dit que deux \'el\'ements $\phi$ et $\psi$ de $\cU$
sont {\it $\Cstar$-\'equivalents}
s'il existe une constante $a \in \Cstar$
telle que $\phi(w) = \psi(aw)$ pour tout $w \in \C$.
On dit aussi que deux \'el\'ements $\phihat$ et $\psihat$ 
de $\widehat{\cU}$
sont {\it $\Cstar$-\'equivalents}
s'il existe une constante $a \in \Cstar$
telle que $\phi_{-n}(w) = \psi_{-n}(aw)$ 
pour tous $n \in \N$ et $w \in \C$.

On note $\cU/\Cstar$ et $\widehat{\cU}/\Cstar$
 les espaces quotients de $\cU$ et $\widehat{\cU}$
par ces $\Cstar$-\'equivalences.
Soit $\cA^{\cK}$ l'adh\'erence de l'ensemble 
$\widehat{\cK}/\Cstar$ dans $\widehat{\cU}/\Cstar$.

En fait, Lyubich et Minsky
ont montr\'e que 
$\widehat{\cU}/\Cstar$ est un feuilletage
par des orbifolds dont 
les rev\^etements universels sont isomorphes \`a $\C$
et que $\cA^{\cK}$ est une lamination 
par des orbifolds du m\^eme type \cite[\S 7]{LM}.

Dans le cas o{\`u} $\cK = \cLM$, on appelle $\cA^{\cLM}$ la {\it lamination de Lyubich-Minsky}
ou {\it LM-lamination} de la fraction rationnelle $f$.
La 3-lamination hyperbolique de Lyubich-Minsky de $f$
est un espace fibr\'e sur $\cA^{\cLM}$.

\parag{Remarque.}
Le choix de $\cK$ est important. 
Par exemple, dans le cas o{\`u} $\cK = \cLM$,
Lyubich et Minsky ont montr\'e que
 $\cA^{\cLM}$ est localement compact
 si $f$ est semi-hyperbolique,
 et ils l'ont utilis\'e pour 
 la d\'emonstration de leur th\'eor\`eme de rigidit\'e.
(La compacit\'e locale n'est pas triviale, car 
l'espace $\cU$ n'est pas localement compact.)

\subsection{Lamination de Zalcman}
Dans la construction de la lamination $\cA^\cK$ ci-dessus, on peut prendre pour $\cK$ l'ensemble des fonctions m\'eromorphes g\'en\'er\'ees par le lemme de Zalcman, 
et ceci permettra de construire une autre lamination, {\it la lamination de Zalcman}.   
En fait, la lamination de Lyubich-Minsky est une sous-lamination de la lamination de Zalcman (\thmref{thm:LM_in_Z}).

\parag{Fonctions de Zalcman.}
Soient $f$ une fonction rationnelle avec $\deg f \ge 2$,
et $z_0 \neq \infty$ un point dans l'ensemble de Julia $J = J(f)$.
D'apr\`es Steinmetz \cite{Ste}, on dit 
qu'un \'el\'ement $\psi \in \cU$  est 
une \textit{fonction de Zalcman} (ou \textit{Z-fonction}) 
de $f$ en $z_0$
si $\psi$ est une limite de la suite  
donn\'ee par le lemme de Zalcman (\lemref{lem:Zalcman}) 
pour la restriction de $\cF = \braces{f^n}_{n \in \N}$ 
au voisinage de $z_0$.
Plus pr\'ecis\'ement, il existe des suites
$\rho_k \to 0, z_k \to z_0, n_k \to \infty$ 
telles que $f^{n_k}(z_k + \rho_k w)$ converge 
vers $\psi(w)$ quand $k \to \infty$ 
uniform\'ement sur tout compact de $\C$.

Soit $\cZ(z_0) = \cZ_f(z_0)$ l'ensemble des  fonctions de Zalcman de $f$ en $z_0$. 
Quand $\infty \in J$, 
on d\'efinit $\cZ(\infty)$ par
$$
\cZ(\infty) := \braces{1/\phi \st \phi \in \cZ_F(0)},
$$
o{\`u} $F$ est la fraction rationnelle d\'efinie par $F(z) := 1/f(1/z)$. 
Alors l'ensemble des fonctions de Zalcman de $f$
est donn\'e par 
$$
\cZ   :=  \bigcup_{z_0\in J} \cZ(z_0).
$$

\parag{Remarque.}
Dans le lemme de Zalcman, les fonctions de la forme 
$\psi_k(w) = f^{n_k}(z_k + \rho_k w)$ 
ne sont pas essentielles pour repr\'esenter
la limite $\psi = \lim \psi_k \in \cZ-\cZ(\infty)$. 
Par exemple, on peut remplacer $\psi_k$ par 
$f^{n_k}(z_k + \rho_k w +e_k(w))$ o{\`u}
 $e_k(w) = o(\rho_k)$ sur tout compact de $\C$. 
Cela sugg\`ere qu'on peut repr\'esenter
toute $\phi \in \cZ$ comme une limite de 
$
f^{n_k} \cc T_k(\rho_k w)
$
o{\`u} $T_k:\Chat \to \Chat$ est 
une rotation sph\'erique telle que $T_k(0) = z_{k}$.

\parag{Invariance de $\cZ$.}
Il n'est pas difficile de montrer:

\begin{prop}\label{prop:Z-func}
L'ensemble $\cZ \subset \cU$ satisfait 
$$
f \cc \cZ=\cZ =\cZ \cc \Aff.
$$ 
\end{prop}
Une d\'emonstration de cette proposition se trouve dans \cite{Ste}.

\parag{Lamination de Zalcman.}
En vertu de la \propref{prop:Z-func},
on peut poser $\cK = \cZ$ 
dans la construction de la lamination $\cA^\cK$ ci-dessus.
On appelle $\cZA$ la \textit{lamination de Zalcman}   
ou \textit{Z-lamination} de $f$.
Alors la Z-lamination contient la LM-lamination:

\begin{thm}\label{thm:LM_in_Z}
L'ensemble $\cLM$ est un sous-ensemble de $\cZ$. 
On a donc l'inclusion $\cA^{\cLM} \subset \cZA$. 
\end{thm}

Pour la d\'emonstration on utilise une caract\'erisation de LM-fonctions:
\begin{prop}[Caract\'erisation de LM-fonctions]\label{prop:LM_LM}
Soient $\zhat = (z_{-n})_{n \ge 0}$ une orbite inverse de $f$
et $\psi$ une fonction m\'eromorphe non constante de $\C$.
Alors les assertions suivantes sont \'equivalentes:
\begin{enumerate}[\rm (1)]
\item
$\zhat \in \cR_\C$ et $\psi \in \cLM(\zhat)$.
\item 
Il existe des suites 
$n_k \in \N$,
 $\rho_k \in \Cstar$ avec $\rho_k \to 0$,
 et des rotations sph\'eriques $T_k:\Chat \to \Chat$ 
 avec $T_k(0) = z_{-n_k}$ telles que
$$
\psi_k(w)  :=  f^{n_k} \cc T_k(\rho_k w)
$$
converge vers $\psi(w)$ sur tout compact de $\C$.
\end{enumerate}
\end{prop}
En vertu du \cite[Lemma 3.6]{KL}, on sait que (1) implique (2). 
En fait, on utilise seulement l'implication 
$(1) \Longrightarrow (2)$ 
dans la d\'emonstration du \thmref{thm:LM_in_Z}.

\begin{pf}[\propref{prop:LM_LM}]
Il suffit de montrer que $(2) \Longrightarrow (1)$. 
Supposons qu'il existe des suites $n_k \to \infty$, $\rho_k \to 0$,
et des rotations sph\'eriques $T_k$ avec $T_k(0)=z_{-n_k}$,
et la fonction $f^{n_k}\cc T_k(\rho_k w)$ 
converge vers $\psi$ sur tout compact de $\C$. 
Comme $\psi(w) = \lim_{k} f^{n_k}\cc T_k(\rho_k w)$ et 
la famille 
$\braces{f^{n_k}\cc T_k(\rho_k w)}_{k \in \N}$
est normale, pour tout $m \in \N$ 
la famille $\braces{f^{n_k-m}\cc T_k(\rho_k w)}_{k \in \N}$ est normale aussi. 
Ainsi on peut trouver une fonction 
$\psi_{-m} \in \cU$ 
telle que $f^m \cc \psi_{-m}= \psi$.
En particulier,
$\psi_{-m}(0) = f^{n_k-m}(z_{-n_k}) = z_{-m}$.

Le degr\'e local $\deg (f^m, z_{-m} )$ est 
uniform\'ement born\'e, car
$$
\deg (\psi, 0 )  =  \deg (f^m, z_{-m} )\deg (\psi_{-m}, 0 ) 
\ge 
\deg (f^m, z_{-m} ).
$$
Donc l'orbite inverse $\zhat$ est r\'eguli\`ere. 

Ensuite on montrera $\zhat \in \cR_\C$:
pour $a \in \C$ posons $\delta_a(w):=  w +a$. 
Alors $(\psi_{-m}\cc \delta_a)_{m \in \N} \in \widehat{\cU}$, 
et l'application $h:a \mapsto (\psi_{-m}\cc \delta_a(0))_{m \in \N}$ 
est une application holomorphe non constante de $\C$ 
sur $L(\zhat)$.
Comme les feuilles de $\cR$ sont isomorphes \`a $\C$, $\D$, 
ou \`a un anneau, $L(\zhat)$ doit \^etre isomorphe \`a  $\C$. 
On a donc $\zhat \in \cR_\C$.

Maintenant on peut prendre une uniformisation $\phi:\C \to L(\zhat)$
avec $\phi(0) = \zhat$. 
Posons $\tilde{\psi}: = \pi_0 \cc \phi \in \cLM(\zhat)$.
En utilisant \cite[Lemma 3.6]{KL}, 
on peut trouver une suite 
$\tilde{\rho}_k \in \Cstar$ telle que 
$\tilde{\psi}(w) = \lim _k \tilde{g}_{k}(w)$ avec 
$\tilde{g}_{k}(w) := f^{n_k}\cc T_k(\tilde{\rho}_k w)$.

Posons ${g}_{k}(w) := f^{n_k}\cc T_k({\rho}_k w)$
et $\lam_k:= \rho_k/\tilde{\rho}_k$
tels qu'on a $g_{k}(w) = \tilde{g}_{k}(\lam_k w)$. 
Comme les suites $g_{k}$ et $\tilde{g}_{k}$ ont leurs limites non constantes $\psi$ et $\tilde{\psi}$, $\lam_k$ converge vers un $\lam \in \Cstar$. 
Ainsi on a $\psi(w) = \tilde{\psi}(\lam w)$ 
o\`u $\tilde{\psi}(\lam w) \in \cLM(\zhat)$,
 et on obtient donc $\psi \in \cLM(\zhat)$.
\QED
\end{pf}

\parag{D\'emonstration du \thmref{thm:LM_in_Z}.}
Soit $\psi$ une LM-fonction. 
Il existe donc une orbite inverse 
$\zhat = (z_{-n})_{n \ge 0} \in \cR_\C$
telle que $\psi \in \cLM(\zhat)$.
Si la suite $\braces{z_{-n}}_{n \ge 0}$ a un point d'adh\'erence 
dans $J = J(f)$, on a clairement $\psi \in \cZ$ par la \propref{prop:LM_LM}, (2).
Sinon tout point d'adh\'erence est 
un point p\'eriodique attractif ou
contenu dans un domaine de rotation.
Dans ces deux cas, toute limite de 
$\psi_k(w)=f^{n_k} \cc T_k(\rho_k w)$ 
est une fonction de $\C$ sur l'ensemble de Fatou. 
Mais c'est contradictoire avec le th\'eor\`eme de Picard.
\QED

\subsection{Co{\"{\i}}ncidence}
Il semble possible que 
l'adh\'erence de $\cLM$ dans $\cU$ et celle de $\cZ$  
co{\"{\i}}ncident, mais on n'a pas de preuve.
On montre ici une condition suffisante pour 
les \'egalit\'es $\cZ = \cLM$ et $\cZA = \cA^{\cLM}$.
En d'autres termes, cette condition permet de 
construire la lamination de Lyubich-Minsky 
\`a partir du lemme de Zalcman. 

\parag{Grande orbite univalente.}
La \textit{grande orbite univalente} ${UGO}(z_0)$ de $z_0$
est l'ensemble des points $\zeta \in \Chat$ tels que 
$f^m(z_0) = f^n(\zeta)$ pour des entiers  $m, n \in \N$ et 
tels qu'il existe une branche $g$ de $f^{-n}\cc f^{m}$ 
d\'efinie dans un voisinage de $z_0$ 
avec $g(z_0) = \zeta$ et $g'(z_0) \neq 0$.

Soit $P$ l'ensemble post-critique de $f$, 
{\it i.e.}, l'adh\'erence des orbites critiques.

\begin{thm}[Co{\"{\i}}ncidence des laminations]\label{thm:Z=LM}
Si $f$ v\'erifie la condition $(\ast)$ suivante, 
on a l'\'egalit\'e $\cZ = \cLM$:
\begin{quote}
$(\ast)$ Pour tout $z_0 \in J$, 
il existe un $z_0' \in UGO(z_0)-P$.
\end{quote}
En particulier, la Z-lamination $\cZA$ co{\"{\i}}ncide
avec la LM-lamination $\cA^{\cLM}$.
\end{thm}
Par exemple, toute fraction rationnelle hyperbolique
satisfait la condition $(\ast)$.
Un peu plus g\'en\'eralement, 
cette condition est satisfaite si la fraction rationnelle est parabolique ({\it i.e.} $J$ n'a pas de point critique).
Une classe d'exemples int\'eressante 
de fonctions satisfaisant $(\ast)$ est 
celle des polyn{\^o}mes quadratiques infiniment renormalisables. 
(Voir \propref{prop_inf_renorm}.)
Donc dans ces cas, on peut construire la LM-lamination 
\`a partir du lemme de Zalcman.

Mais la condition $(\ast)$ n'est pas v\'erifi\'ee 
si $f$ est un polyn{\^o}me de Chebychev ou 
un exemple de Latt\`es. 
(Par exemple, posons $f(z) = z^2-2$ et $z_0 = 2$.)
Si on a $J = P$, alors $f$ ne v\'erifie pas $(\ast)$ aussi.

Pour d\'emontrer le th\'eor\`eme, 
on montre:

\begin{lem}\label{lem_ZLM1}
Pour tous $z_0 \in J$ et $\zeta_0 \in UGO(z_0)$, 
on a $\cZ(z_0) = \cZ(\zeta_0)$. 
\end{lem}

\begin{pf}
En utilisant une conjugaison,
nous pouvons supposer que 
$z_0$ et $\zeta_0$ ne sont pas $\infty$.
Prenons une fonction $\psi(w) = \lim_{k \to \infty} f^{n_k}(z_k + \rho_k w) \in \cZ(z_0)$ avec 
$n_k \to \infty$, $z_k \to z_0$, et $\rho_k \to 0$.
Soit $g$ la branche univalente de $f^{-n} \cc f^m$
telle que $g(z_0) = \zeta_0$. 
Par un th\'eor\`eme de Koebe, 
on a 
$g(z) = \zeta_k + A_k (z-z_k) + o(z-z_k)$ 
o{\`u} $\zeta_k = g(z_k)$ et $A_k = g'(z_k) \asymp g'(z_0) \neq 0$ 
pour tout $k \gg 0$. 
Par suite
$$
f^{n_k} \cc g^{-1} \cc g(z_k + \rho_k w)
 =  f^{n_k-m + n}(\zeta_k + A_k \rho_k w + o(\rho_k)) 
$$
dans tout compact du $w$-plan. 
Comme $\zeta_k \to \zeta_0$ quand $k \to \infty$, on obtient $\psi \in \cZ(\zeta_0)$.

\QED\end{pf}

Ensuite on montre:

\begin{lem}\label{lem_ZLM2}
Pour tout $\zeta_0 \in J-P$, 
l'ensemble $\cZ(\zeta_0)$ est un sous-ensemble de $\cLM$. 
\end{lem}

\begin{pf}
Soit $\psi(w):= \lim f^{n_k}(\zeta_k + \rho_k w) \in \cZ(\zeta_0)$
avec $n_k \to \infty$, $\zeta_k \to \zeta_0$, et $\rho_k \to 0$.
Supposons d'abord qu'on puisse extraire 
une sous-suite telle que $\zeta_k \in J$ pour tout $k \in \N$.  
L'ensemble de Julia $J$ peut \^etre 
arbitrairement approch\'e par $f^{-m}(\{z\})$ 
dans la topologie de Hausdorff
quand $z \in J$ et $m \to \infty$.
Donc on peut choisir un point $z_0 \in J$ et une orbite inverse 
$\zhat = (z_{-m})$ tels que 
pour tout $k$ il existe un $m_k \in \N$
suffisamment grand par rapport  \`a $n_k$ 
tel que $z_{-n_k} = f^{m_k-n_k}(z_{-m_k}) = \zeta_k + o(|\rho_k|)$. 
(Car $J$ n'a pas de point isol\'e.) 

Comme $\zeta_0 \notin P$ par hypoth\`ese, 
il existe un disque $D_0$ centr\'e en $\zeta_0$
et une branche univalente $g_k:D_0 \to \Chat$ de $f^{-(m_k-n_k)}$
telle que $g_k(z_{-n_k}) = z_{-m_k}$. 
Le multiplicateur $\lam_k:= g_k'(\zeta_0)$ satisfait 
$\liminf \lam_k = 0$ par la normalit\'e de $\{g_k|_{D_0}\}$. 
Fixons un $r >0$ arbitrairement grand.
Pour tout $k \gg 0$ (assez grand),
le disque $D_k: = \D(\zeta_k, r|\rho_k|)$ 
est contenu dans $D_0$ et $z_{-n_k} \in D_k$.
Par un th\'eor\`eme de Koebe,
pour tout $z \in g_k(D_k)$, on a $|z-z_{-m_k}| = O(\rho_k\lam_k)$. 
Donc on peut \'ecrire $f^{m_k-n_k}|_{g_k(D_k)}$ comme
$Z = f^{m_k-n_k}(z) = z_{-n_k}  + a_k(z-z_{-m_k}) + O(\rho_k^2\lam_k^2)$ avec $a_k \asymp \lam_k^{-1} \to \infty$ (en utilisant 
le th\'eor\`eme de Koebe encore une fois). 
Par suite, on a 
$g_k(Z) = z_{-m_k} + (Z-\zeta_k)/a_k + o(|\rho_k/a_k|)$ et
$$
\psi(w)  =  \lim f^{m_k}\cc g_k(\zeta_k + \rho_k w) 
 =  \lim f^{m_k}(z_{-m_k} + \rho_k/a_k (w + o(1))).
$$
On obtient donc $\psi \in \cLM$. 

Ensuite, supposons que $\zeta_k \notin J$ 
pour tout $k \in \N$ assez grand. 
Soit $R_k: = \dist_{\C}(\zeta_k,J)$.
Si $R_k > r |\rho_k|$ pour tous $r>0$ et $k$ assez grand,
l'ensemble de Fatou contient le disque 
$D_k =  \D(\zeta_k, r|\rho_k|)$, 
et donc $\psi$ 
est une fonction de $\C$ sur l'ensemble de Fatou. 
C'est contradictoire avec le th\'eor\`eme de Picard.
Donc il existe un $r_0>0$ et un nombre infini de $k$
tels que $R_k \le r_0|\rho_k|  \to 0~(k \to \infty)$. 
Prenons $\zeta_k' \in J$ avec 
$|\zeta_k-\zeta_k'| = R_k$.
Soit $\delta_k: = \zeta_k-\zeta_k'$.
Comme $|\delta_k/\rho_k| =R_k/|\rho_k| \le r_0$,
on peut trouver 
un $a \in \overline{\D(r_0)}$ et une sous-suite 
de $\braces{k}$ 
avec $\delta_k/\rho_k \to a$ quand $k \to \infty$.
Alors
$$
f^{n_k}(\zeta'_k + \rho_k(w + a))
 = f^{n_k}(\zeta'_k + \rho_k(w + \delta_k/\rho_k + o(1)))
 = f^{n_k}(\zeta_k + \rho_kw + o(\rho_k))
$$
et on a donc $\lim_{k\to\infty} f^{n_k}(\zeta'_k + \rho_kw) = \psi(w-a)$. 
Comme $\zeta_k' \in J$ et $\zeta_k' \to \zeta_0$,
on obtient $\psi(w-a) \in \cLM$, et donc $\psi(w) \in \cLM$.
(Rappelons que $\cLM = \cLM \cc \Aff$.)
\QED\end{pf}

\parag{D\'emonstration du \thmref{thm:Z=LM}.}
Comme $\cLM \subset \cZ$ (le \thmref{thm:LM_in_Z}), 
il suffit de montrer  
que $\cZ(z_0) \subset \cLM$ pour tout $z_0 \in J$.
Supposons que la fonction $f$ v\'erifie la condition $(\ast)$ et fixons $z_0 \in J$ avec $z_0' \in UGO(z_0)-P$. 
Alors on obtient le th\'eor\`eme en posant $\zeta_0:=z_0'$
dans ces deux lemmes. 
\QED

\parag{Question.}
Le probl\`eme suivant n'est pas encore r\'esolu:  
\textit{
Sous quelle condition n\'ecessaire et suffisant  
a-t-on l'\'egalit\'e $\cZ = \cLM$ (ou $\cZA = \cA^{\cLM}$) ?
}

Remarquons que, en fait, la condition $\zeta_0 \in UGO(z_0)$ du \lemref{lem_ZLM1} n'est pas optimale pour avoir $\cZ(z_0) = \cZ(\zeta_0)$. Voir \cite[Theorem 4]{Ste}. Donc on peut remplacer $(\ast)$
par une condition un peu plus faible (mais plus compliqu\'ee !).

\subsection{Polyn{\^o}mes quadratiques infiniment renormalisables}
On donne une d\'emonstration de:

\begin{prop}\label{prop_inf_renorm}
Les polyn\^omes quadratiques $f_c(z) = z^2 + c$ 
infiniment renormalisables
satisfont la condition $(\ast)$.
Donc on a $\cA^{\cLM} = \cA^{\cZ}$.
\end{prop}

\parag{D\'emonstration.}
Supposons que $f(z) = z^2 + c$ est infiniment renormalisable.
Comme tout les points p\'eriodiques de $f$ sont r\'epulsifs 
et l'ensemble post-ctirique $P$ ne contient pas 
les points p\'eriodiques
 (voir \cite[Thm. 8.1]{Mc1}), 
il existe un voisinage $U$ du point fix\'e r\'epulsif $\beta \in J$ 
de $f$ avec $U \cap P = \emptyset$. 
Donc il existe un $n \in \N$ tel que $J \subset f^n(U)$. 

Pour $z_0 \in J$, s'il existe un $w_0 \in U$ tel que 
$f^n(w_0) = z_0$ et $(f^n)'(w_0) \neq 0$, 
alors $w_0 \in UGO(z_0)-P$. 
Sinon $f^n(w_0) = z_0$ et $(f^n)'(w_0) = 0$. 
Dans ce cas il y a un entier $k \le n$ tel que $z_0$ est 
l'image du point critique $z = 0$ par $f^k$. 
On peut alors trouver une renormalisation 
$g = f^p|_{V'}:V' \to V$ avec $z_0 \in V' \Subset V$, 
dont le petit ensemble de Julia $J_0 = \bigcap_{j \ge 0}g^{-j}(V) \subsetneq J$ satisfait $f^{p}(J_0) = J_0$; $z_0 \in J_0$; et $P \subset \bigcup_{i = 1}^p f^i(J_0)$. (Voir par exemple \cite[\S 10.1]{Mc1}.)
Alors $J_0$ ou $f(J_0)$ ne contient pas 
la valeur critique $f(0) = c$, 
et il existe un des composantes $L$ de $f^{-1}(J_0)$ ou $f^{-1}(f(J_0))$ tel que $L \neq f^i(J_0)~(1 \le i \le p)$ et $L \cap P =\emptyset$. 
Donc on peut trouver un $z'_0 \in L$ avec $f(z_0') = z_0$ ou $f(z'_0) =f(z_0)$, et $z_0' \in UGO(z_0)-P$. 
\QED

\subsection{Une question de Steinmetz}

Dans \cite{Ste}, Steinmetz pose plusieurs questions 
sur l'espace $\cZ$ des Z-fonctions.
Une de ses questions est la suivante:
\textit{
Pour une fonction $f$ donn\'ee,
peut-on d\'eterminer l'ensemble 
$\braces{\zeta \in J \st \cZ = \cZ(\zeta)}$  ?
}
Comme un corollaire de la d\'emonstration du \thmref{thm:Z=LM},
on donne une condition suffisante 
pour que l'ensemble soit $J$:  

\begin{thm}\label{thm:St-Q}
Si $f$ satisfait $(\ast)$ du \thmref{thm:Z=LM}, 
alors on a $\cZ = \cZ(\zeta)$ pour tout $\zeta \in J$.
\end{thm}

C'est une am\'elioration du \cite[Thm.5]{Ste}.

\begin{pf}
Par le \lemref{lem_ZLM1}, pour tout $z_0 \in J$,
il existe un $\zeta_0 \in (J-P) \cap UGO(z_0)$ 
avec $\cZ(z_0) = \cZ(\zeta_0)$. 

Fixons un point arbitraire $\zeta$ dans l'ensemble de Julia.
Alors on peut prendre l'orbite inverse $\zhat = (z_{-n})$ 
comme dans la d\'emonstration du \lemref{lem_ZLM2}
telle que $z_{-m_k}$ converge vers $\zeta$.
Par suite $\cZ(z_0) = \cZ(\zeta_0) \subset \cZ(\zeta)$. 
En prenant la r\'eunion sur tous les points $z_0 \in J$, 
on obtient $\cZ \subset \cZ(\zeta)$. 
L'autre inclusion est \'evidente.
\QED
\end{pf}

\section{Caract\'erisations de points coniques}
Ensuite on consid\`ere 
des points coniques dans l'ensemble de Julia 
pour une fraction rationnelle $f$ de degr\'e au moins deux.  
Il y a plusieurs d\'efinitions des points coniques
(voir \cite{Pr}),
et chacune d'elles est li\'ee \`a 
la propri\'et\'e d'expansion dans l'ensemble de Julia.
Dans cette section, 
on compare la notion de point conique
d\'efinie par Lyubich et Minsky \cite{LM}, 
avec celle par Martin et Mayer \cite{MM} 
dans le contexte unifi\'e des laminations.

\subsection{Points LM-coniques}

\parag{Dynamique sur la lamination.}
Les points coniques de Lyubich-Minsky sont d\'efinis
en terme de la dynamique sur la LM-lamination $\cA^{\cLM}$.
Il existe une action inversible  
$\fhat:\widehat{\cU} \to \widehat{\cU}$ 
d\'efinie par
$\fhat:(\psi_{-n}) \mapsto (f \cc \psi_{-n})$. 
Les LM- et Z-laminations h\'eritent de cette action,
not\'ee $\fhat:\cA^{\cLM} \to \cA^{\cLM}$ et $\fhat:\cZA \to \cZA$.
(En fait, la projection canonique de $\widehat{\cU}$ 
sur $\widehat{\cU}/\Cstar$ et la post-composition par $f$ sont commutatifs 
pour tout \'el\'ement de $\widehat{\cU}$.)

Pour tout $[\psihat] \in \widehat{\cU}/\Cstar$ 
avec un repr\'esentant 
$\psihat =(\psi_{0}, \psi_{-1}, \ldots) \in \widehat{\cU}$, 
la valeur $\psi_{0}(0)$ dans $\Chat$ 
ne d\'epend pas du choix du repr\'esentant. 
On dit que l'application $\hat{\pi}: \widehat{\cU}/\Cstar \to \Chat$ 
d\'efinie par $\hat{\pi}([\psi]) := \psi_{0}(0)$ est
la \textit{projection}.
Remarquons que la dynamique 
$\fhat:\widehat{\cU}/\Cstar \to \widehat{\cU}/\Cstar$ 
est semi-conjugu\'ee \`a la dynamique originale 
$f:\Chat \to \Chat$ par $\hat{\pi}$.

\parag{Points LM-coniques.}
En analogie avec les groupes kleiniens,
Lyubich et Minsky ont d\'efini les \textit{points coniques} 
dans $\cA^{\cLM} \cap \hat{\pi}^{-1}(J)$, o\`u $J = J(f)$: 

\begin{df}[Points coniques dans la lamination, {\cite[\S 8]{LM}}]
Un point $\tilde{z}$ dans $\cA^{\cLM} \cap \hat{\pi}^{-1}(J)$
est appel\'e \textit{conique} si l'orbite 
$\{\tilde{z}, \fhat(\tilde{z}), \fhat^2 (\tilde{z}), \ldots \}$ 
a un point d'adh\'erence dans $\cA^{\cLM}$.
\end{df}

Dans cet article, nous appliquons cette d\'efinition 
\`a la Z-lamination $\cZA$ \`a la place de $\cA^{\cLM}$.
(Comme elles co\"{\i}ncident souvent par \thmref{thm:Z=LM}.)
De plus, on d\'efinit les points coniques dans l'ensemble de Julia:

\begin{df}[Points coniques de Lyubich-Minsky, modifi\'es]
Un point $z_0$ dans l'ensemble de Julia 
est appel\'e \textit{LM-conique} 
s'il existe un \'el\'ement $\tilde{z} \in \cZA$ 
tel que $\hat{\pi}(\tilde{z}) = z_0$,
et tel que l'orbite 
$\{\tilde{z}, \fhat(\tilde{z}), \fhat^2 (\tilde{z}), \ldots \}$ 
 ait un point d'adh\'erence dans $\cZA$.
\end{df} 

\parag{}
Pour distinguer d'autres d\'efinitions de point conique,
on dit que le point conique ci-dessus est \textit{LM-conique}, 
et on note $\Lam_{LM} \subset J$ l'ensemble des points LM-coniques. 

Une propri\'et\'e importante de $\Lam_{LM}$ est:

\begin{thm}[Lyubich-Minsky, rigidit\'e de $\Lam_{LM}$]
\label{thm:LM}
Si $\Lam_{LM}$ d'une fraction rationnelle $f$ 
a un champ de droites invariant, 
alors $f$ est un exemple de Latt\`es.
\end{thm}
Pour les terminologies, voir par exemple \cite{Ha, MM}. 
Pour la d\'emonstration, voir \cite[Prop. 8.9]{LM}.
(Malgr\'e l'usage de $\cZA$ \`a la place de $\cA^{\cLM}$,
la d\'emonstration est la m\^eme.) 
Lyubich et Minsky ont aussi montr\'e
que toute fraction rationnelle semi-hyperbolique 
(\textit{convexe co-compacte} dans leur terminologie)
satisfait $J = \Lam_{LM}$ et on a la rigidit\'e sur l'ensemble de Julia.

\parag{Autres caract\'erisations des points LM-coniques.}
On donnera d'autres caract\'erisations importantes de $\Lam_{LM}$:

\begin{thm}[Points LM-coniques]\label{thm_another_LM}
Pour $z_0 \in J$, les assertions suivantes sont \'equivalentes:
\begin{itemize}
\item[\rm (LM1)]
il existe une fonction $\phi \in \cZ(z_0)$ de la forme 
$$
\phi(w)  =  \lim_{k \to \infty}f^{n_k}(z_0 + \rho_k w)
$$
o{\`u} la convergence est uniforme sur tout compact de $\C$.
\item[\rm (LM2)]
$z_0$ est LM-conique. 
\item[\rm (LM3)]
pour tout \'el\'ement $\tilde{z} \in \cZA$ avec $\hat{\pi}(\tilde{z}) = z_0$, 
l'orbite 
$\{\tilde{z}, \fhat(\tilde{z}), \fhat^2 (\tilde{z}), \ldots \}$ 
a un point d'adh\'erence dans $\cZA$.
\end{itemize}
\end{thm}
Si $z_0 = \infty$ dans (LM1), 
on utilise $F(z) = 1/f(1/z)$ au lieu de $f$.

\begin{pf}
D'abord on montre (LM2) $\Longrightarrow$ (LM1):
supposons que $\tilde{z} \in \cZA$ satisfasse $\fhat^{n_k}(\tilde{z}) \to \tilde{\zeta} \in \cZA$. 
Soient $\hat{\psi} = (\psi_{-n})$ et $\hat{\phi} = (\phi_{-n})$ 
des repr\'esentants de $\tilde{z}$ et $\tilde{\zeta}$ 
dans $\widehat{\cU}$. 
Alors il existe $\lam_k \in \Cstar$ 
tel que pour tout $N \in \N$ fix\'e,
$$
f^{n_k} \cc \psi_{-N}\cc \lam_k ~\to~ \phi_{-N}~~~(k \to \infty)
$$
uniform\'ement sur tout disque $D \subset \C$.
(Ici $\lam_k$ est l'application $\lam_k: w \mapsto \lam_k w$.)
Si $\lam_k \asymp 1$,
il existe un ouvert $U$ tel que 
$U \cap J \neq \emptyset$ et 
$U \subset \psi_{-N} \cc \lam_k (D)$ 
pour tout $k \gg 0$ (o\`u $N$ est fix\'e),
et donc on n'a pas la convergence ci-dessus. 
Par suite $\lam_k$ tend vers $0$. 

Posons $N = 0$. 
Alors il existe un $A_0 \neq 0$ et $p \ge 1$ tels que
$
\psi_0(z) = z_0 + A_0 z^p +o(z^p)
$
pr\`es $z=0$ et on a donc 
$$
\psi_0(\lam_k w)  =  z_0 + A_0  \lam_k^p w^p +o(\lam_k^p)
$$
dans $D$ pour tout $k \gg 0$. Du coup
$$
\phi_0(w)  =  
\lim_{k \to \infty}
f^{n_k}\cc \psi_0 \cc \lam_k (w) 
 = 
\lim_{k \to \infty} f^{n_k}(z_0 + A_0  \lam_k^p w^p).
$$
Donc il existe  
$\phi(w) = \lim_{k \to \infty} f^{n_k}(z_0 + {\rho}_k w) \in \cZ(z_0)$ avec $\rho_k = A_0  \lam_k^p \to 0$ et $\phi_0(w) = {\phi}(w^p)$. 

Ensuite on montre (LM1) $\Longrightarrow$ (LM3).
Prenons $\tilde{z} \in \cZA$ avec $\hat{\pi}(\tilde{z}) = z_0$ 
et un de ses repr\'esentants $\hat{\psi} = (\psi_{-n}) \in \widehat{\cU}$. 
Alors on a  
$$
\psi_0(w) = z_0 + A_0 w^p +o(w^p)
$$
pr\`es de $w = 0$. 
On peut supposer que $A_0 = 1$ en rempla{\c{c}}ant
le repr\'esentant $\hat{\psi}$ par 
$\hat{\psi} \cc A_0^{1/p}$.
Par hypoth\`ese, 
il existe une $\phi(w) = \lim_k f^{n_k}(z_0+\rho_k w)\in \cZ$ o\`u $n_k \to \infty$, $\rho_k \to 0$.
On a donc
$$
\lim_{k \to \infty} f^{n_k}\cc \psi_0 \cc \rho_k^{1/p} (w)
 =  
\lim_{k \to \infty} f^{n_k}(z_0 + \rho_k w^p + o(\rho_k))
 = 
\phi(w^p) ~ =: \phi_0(w)
$$ 
dans tout grand disque.
En d'autres termes, 
la suite $\{f^{n_k}\cc \psi_0 \cc \rho_k^{1/p}\}_{k \in \N} \subset \cZ$ converges vers $\phi_0(w) \in \cU$.
Comme l'espace $\cZA$ est ferm\'e,
la suite $[f^{n_k}\cc \psi_0] \in \cZA$ 
converge vers $[\phi_0] \in \cZA$ 
dans la topologie de $\cU/\Cstar$.
On peut appliquer le m\^eme argument \`a tout $\psi_{-n}$,
et ainsi on peut trouver une limite de $\fhat^{n_k}(\tilde{z})$ 
dans $\cZA$.

Enfin, l'implication (LM3) $\Longrightarrow$ (LM2) est \'evidente. 
\QED
\end{pf}

\parag{Remarque.}
D'apr\`es la d\'emonstration tout point d'adh\'erence  
de l'orbite $\{\fhat^n (\tilde{z})\}_{n \in \N}$ dans (LM3) 
a un repr\'esentant de la forme 
$\phi_{-m}(w^p) = \lim_{k \to \infty} f^{n_k-m}(z_0 + \rho_k w^p)$. 
Comme Lyubich et Minsky ont remarqu\'e,
cela est consid\'er\'e comme un effet
de la vari\'et\'e stable locale de la dynamique $\fhat:\cZ \to \cZ$.

\subsection{Points MM-coniques}

\parag{Th\'eor\`eme de Ha{\"{\i}}ssinsky-Martin-Mayer sur la rigidit\'e.}
Il y a un am\'elioration du th\'eor\`eme de Lyubich-Minsky
d{\^u} \`a Ha\"{\i}ssinsky \cite{Ha}, Martin et Mayer \cite{MM}. 
D'abord on pr\'esente une autre d\'efinition 
du point conique dans \cite{MM}.
(L'importance de la condition ci-dessous 
\'etait aussi remarqu\'ee par 
Astala-Ha{\"{\i}}ssinsky \cite{Ha} et Steinmetz \cite{Ste}.)

\begin{df}[Points coniques de Martin-Mayer]
On dit qu'un point $z_0 \neq \infty$ dans l'ensemble de Julia
est \textit{MM-conique} s'il existe des suites $n_k \in \N$ et $\rho_k \in \Cstar$ avec $\rho_k \to 0$ telles que 
la fonction $\psi_k(w) = f^{n_k}(z_0 + \rho_k w)$ converge vers
une fonction m\'eromorphe non constante $\psi:\D \to \Chat$ 
uniform\'ement sur tout compact de $\D$.
\end{df} 
On dit aussi que $z_0 = \infty$ est MM-conique si 
$0$ est un point MM-conique de la fraction rationnelle 
$F(w)=1/f(1/w)$. (Alors la notion de point MM-conique
est invariante par des changements de coordonn\'ees.)
On note $\Lam_{MM}$ l'ensemble des points MM-coniques. 

Remarquons que cette d\'efinition 
est une variante de (LM1). 
La diff\'erence essentielle est que
le domaine de $\psi = \lim \psi_k$ ne doit pas \^etre $\C$.

\begin{prop}[Caract\'erisation topologique des points MM-coniques]
\label{prop:MM-conical}
Pour $z_0 \in J$, les assertions suivantes sont \'equivalentes:
\begin{itemize}
\item[\rm (MM0)] 
il existe un $d \in \N$ et un $r >0$ tels que 
pour $D_n:=\B(f^n(z),r)$ et la composante connexe $D_n'$ de 
$f^{-n}(D_n)$ contenant $z_0$, 
$\deg(f^n:D_n' \to D_n) \le d$ pour un nombre infini de $n \in \N$.   \item[\rm (MM1)]
$z_0$ est MM-conique.
\end{itemize}
\end{prop}
On note ici $\B(x,r)$ le disque sph\'erique de rayon $r$ dont le centre est $x$.  
Pour la d\'emonstration, voir \cite{Ha, MM}.
Comme (MM1) est plus faible que (LM1), on a  
\begin{prop}\label{prop:inclution}
$\Lam_{LM} \subset \Lam_{MM}$.
\end{prop}
L'implication (LM2) $\Longrightarrow$ (MM0) 
est \'equivalente \`a la proposition,
et d\'ej\`a remarqu\'e dans \cite[Prop.8.7]{LM}.

Voici une am\'elioration du \thmref{thm:LM} (\cite[Prop. 5.2]{Ha} and \cite[Thm.1.2]{MM}):

\begin{thm}[Rigidit\'e par Ha\"{\i}ssinsky, Martin-Mayer]\label{thm:HMM}
Si $\Lam_{MM}$ d'une fraction rationnelle $f$ 
a un champ de droites invariant, 
alors $f$ est un exemple de Latt\`es.
\end{thm}

Dans la suite, pour comparer les d\'efinitions 
de point conique
on introduit un nouvel espace topologique
dans l'esprit de la lamination de Lyubich-Minsky.

\subsection{Espace des germes et points MM-coniques}

\parag{L'espace des germes.}
Nous essayons de modifier la construction de la $Z$-lamination.
L'id\'ee essentielle est 
\og la topologie des germes \fg\ 
introduite par \cite{KL}. 
(On peut aussi trouver une id\'ee similaire dans \cite{BS}.)  
Soit $\cV$ l'ensemble des fonctions m\'eromorphes
non constantes
d\'efinies (au moins) dans un disque centr\'e \`a l'origine.
Nous consid\'erons tout \'el\'ement $\psi \in \cV$ comme un germe en $0$, donc on ignore son domaine maximal.

On dit qu'une suite $\psi_n \in \cV$ converge vers $\psi \in \cV$ 
s'il existe un disque $D$ centr\'e \`a $w=0$ 
tel que $\psi|_D$ et $\psi_n|_D$ soient d\'efinies pour tout $n \gg 0$, et que $\psi_n$ converge vers $\psi$ uniform\'ement dans $D$.

Soit $f$ une fraction rationnelle.
{\'E}videmment $f \cc \psi(w) \in \cV$ et $\psi(aw) \in \cV$ 
pour tous $a \in \Cstar$ et $\psi \in \cV$. 
On a donc une action $f \cc:\cV \to \cV$,
 et le quotient $\cV/\Cstar$ est bien d\'efini.

Nous consid\'erons maintenant $\cZ$ 
comme un sous-ensemble de $\cV$.
On d\'efinit \textit{l'espace des germes de Zalcman}
par l'adh\'erence de $\cZ/\Cstar$ dans $\cV/\Cstar$,
not\'e par $\cZG$. 
C'est-\`a-dire que $\cZG$ est l'ensemble de 
tous les germes g\'en\'er\'es par 
le principe du lemme de Zalcman.
(On peut construire l'espace des germes dans 
$\cV^\N/\Cstar$ avec une dynamique inversible, 
mais la caract\'erisation de points MM-coniques dans 
le \thmref{thm_another_MM} sera presque la m\^eme.)

Un m\'erite de cet espace est d'\'etendre 
la condition de la convergence dans $\cU$. 
En fait, dans $\cU$ la convergence est celle sur le plan $\C$, 
mais dans $\cV$ un disque $\D$ est suffisant.
De plus, son rayon est changeable,
car on prend le quotient par l'action de $\Cstar$.

\parag{Caract\'erisation de points MM-coniques.}
Enfin, on montre:

\begin{thm}[L'analogue du \thmref{thm_another_LM}]
\label{thm_another_MM}
Pour $z_0 \in J$ les assertions suivantes sont \'equivalentes:
\begin{itemize}
\item[\rm (MM1)]
$z_0$ est MM-conique.
\item[\rm (MM2)]
Il existe un \'el\'ement $\psi \in \cV$ 
avec $\psi(0)=z_0$ tel que pour $\tilde{z}:=[\psi] \in \cZG$, 
l'orbite $\braces{\tilde{z}, f(\tilde{z}), f^2 (\tilde{z}), \ldots }$ 
ait un point d'adh\'erence dans $\cZG$. 
\item[\rm (MM3)]
Pour tout $\tilde{z}=[\psi] \in \cZG$ avec $\psi(0)=z_0$, 
l'orbite $\braces{\tilde{z}, f(\tilde{z}), f^2 (\tilde{z}), \ldots }$ 
a un point d'adh\'erence dans $\cZG$. 
\end{itemize}
\end{thm}

\begin{pf}
La d\'emonstration est parfaitement analogue 
\`a celle du \thmref{thm_another_LM}. 
On doit seulement remplacer $\C$, le domaine de d\'efinition
des fonctions, par un disque centr\'e  \`a l'origine.
\QED
\end{pf}

\parag{Remarque.}
On peut donner une d\'emonstration du \thmref{thm:HMM} 
en suivant le m\^eme argument que celui du \thmref{thm:LM} (\cite[Proposition 8.9]{LM}) dans le contexte des germes.

\end{document}